\documentclass[preprint,1p,leqno]{elsarticle}
\usepackage{amsmath}
\usepackage{amssymb}
\usepackage{mathrsfs}
\usepackage{color}
\usepackage{marvosym}

\usepackage[]{graphicx}
\usepackage{psfrag}
\usepackage{pdfpages}
\usepackage[dvipsnames]{xcolor}

\usepackage[makeroom]{cancel}

\usepackage{amsthm}

\usepackage{bigstrut}
\usepackage{multirow}
\usepackage{subfig}
\usepackage{easyReview}
\setreviewson
\usepackage{hyperref}

\usepackage{geometry}
\newcommand{\bu}{\boldsymbol{u}}
\newcommand{\bx}{\xi}
\newcommand{\bn}{\boldsymbol{n}}
\newcommand{\bF}{\boldsymbol{f}}

\newcommand{\bv}{\boldsymbol{v}}

\newcommand{\Rn}{\mathrm{Re}}

\journal{arXiv}
\begin{document}
\begin{frontmatter}
\title{A MEEVC discretization for two-dimensional incompressible Navier-Stokes equations with general boundary conditions}

\author[1]{Yi Zhang\corref{cor1}}
\ead{zhangyi_aero@hotmail.com}
\author[2]{Artur Palha}
\ead{a.palha@tudelft.nl}
\author[2]{Marc Gerritsma}
\ead{m.i.gerritsma@tudelft.nl}
\author[3]{Qinghe Yao}
\ead{yaoqhe@sysu.edu.cn}
%\author[1]{Stefano Stramigioli}
%\ead{s.stramigioli@utwente.nl}

\affiliation[1]{organization={University of Twente},
	addressline={Drienerlolaan 5},
	postcode={7522 NB},
	city={Enschede},
	country={the Netherlands}}

\affiliation[2]{organization={Delft University of Technology},
addressline={Mekelweg 5},
postcode={2628 CD},
city={Delft},
country={the Netherlands}}

\affiliation[3]{organization={Sun Yat-sen University},
	addressline={Gongchang road 66},
	postcode={518107},
	city={Shenzhen},
	country={China}}
\cortext[cor1]{Corresponding author}

%\fntext[fn1]{This is the first author footnote.}

\begin{abstract}
	In this work, we introduce a mass, energy, enstrophy and vorticity conserving (MEEVC) mixed finite element discretization for two-dimensional incompressible Navier-Stokes equations as an alternative to the original MEEVC scheme proposed in [A. Palha and M. Gerritsma, J. Comput. Phys., 2017].
%	[\emph{A. Palha and M. Gerritsma, A mass, energy, enstrophy and vorticity conserving (MEEVC) mimetic spectral element discretization for the 2D incompressible Navier-Stokes equations, 2017}]. 
	The present method can incorporate no-slip boundary conditions.
	Conservation properties are proven. Supportive numerical experiments with both exact and inexact quadrature are provided. % directly without \Yi{compromising on} the vorticity conservation.
\end{abstract}
\end{frontmatter}

%_____________________________________________________________________________________________________
%>>>>>>>>>>>>>>>>>>>>>>>>>>>>>>>>>>>>>> SECTION <<<<<<<<<<<<<<<<<<<<<<<<<<<<<<<<<<<<<<<<<<<<<<<<<<<<<<<
\section{Introduction}
We consider the dimensionless rotational or Lamb form of two-dimensional incompressible (or, more strictly speaking, constant density) Navier-Stokes equations in a space-time domain, see for example \cite{zang1991rotation, layton2009accuracy, ZHANG2002764},
\begin{subequations}\label{Eq: NSE}
\begin{align}
\partial_{t}\bu +  \omega\times\bu +\Rn^{-1}\nabla\times\omega + \nabla P &= \bF &&\text{in}\ \Omega\times(0,T] ,\label{Eq: NSE a}\\
\omega-\nabla\times\bu&=\boldsymbol{0} &&\text{in}\ \Omega\times(0,T],\label{Eq: NSE b}\\
\nabla\cdot\bu &= 0&&\text{in}\ \Omega\times(0,T],\label{Eq: NSE c}
\end{align}
\end{subequations}
where $ \Omega\subset\mathbb{R}^2 $, is a simply connected, bounded domain with a Lipschitz boundary $ \partial\Omega $, $ \partial_{t} :=\frac{\partial}{\partial t}$, $ \bu $ is the velocity field, $ \omega $ is the vorticity field, $ P:=p + \frac{1}{2}\bu\cdot\bu $ is the total pressure (with $ p $ being the static pressure), $ \bF $ is the external (body) force, and $ \Rn $ denotes the Reynolds number. \eqref{Eq: NSE} is supplemented with an initial condition, 
\begin{equation}\label{Eq initial conditon}
	\bu^{0} = \bu(\boldsymbol{x},t_0) ,
\end{equation}
and two pairs of boundary conditions,
\begin{equation}\label{Eq: BCs}
	\begin{aligned}
		\left\lbrace
		\begin{aligned}
			\bu\cdot\bn &= \widehat{u}_{\perp}\qquad   &&\text{on}\ \Gamma_{\perp}\times(0,T] \\
			P &= \widehat{P} \qquad&&\text{on}\ \Gamma_{\widehat{P}}\times(0,T] 
		\end{aligned}\,,
	\right.\\
	\left\lbrace
	\begin{aligned}
		\omega  &= \widehat{\omega} \qquad&&\text{on}\ \Gamma_{\widehat{\omega}}\times(0,T] \\
		\bu\times\bn &= \widehat{u}_{\parallel} \qquad&&\text{on}\ \Gamma_{\parallel}\times(0,T] 
	\end{aligned}\,,
\right.
\end{aligned}
\end{equation}
where $ \bn $ is the unit outward norm vector. In each pair, the boundary sections are disjoint, and therefore the boundary conditions are not active at the same boundary section, i.e., $  \Gamma_{\perp}\cap \Gamma_{\widehat{P}} =  \Gamma_{\widehat{\omega}}\cap \Gamma_{\parallel} = \emptyset $, and each pair is active over the whole boundary, $ \Gamma_{\perp}\cup \Gamma_{\widehat{P}} = \Gamma_{\widehat{\omega}}\cup \Gamma_{\parallel} = \partial\Omega $, also see \cite[Table~1]{kreeftstokes}. Note that, when $ \Gamma_{\widehat{P}} =\emptyset $, there is a singular mode in the total pressure of \eqref{Eq: NSE}; the total pressure is determined up to a constant. 

The divergence free condition of velocity, \eqref{Eq: NSE c}, implies mass conservation. For two-dimensional incompressible flows, integral quantities of interest include
\[
\begin{aligned}
	\emph{total kinetic energy (or simply energy) }&\mathcal{K}:=\dfrac{1}{2}\int_{\Omega} \bu\cdot\bu\ \mathrm{d}\Omega, \\
	\emph{total enstrophy (or simply enstrophy) }&\mathcal{E}:=\dfrac{1}{2}\int_{\Omega}\omega\cdot\omega\ \mathrm{d}\Omega,\\
	\emph{total palinstrophy (or simply palinstrophy) } & \mathcal{P} := \dfrac{1}{2}\int_{\Omega} \left( \nabla \times \omega\right)  \cdot \left( \nabla \times \omega\right)  \ \mathrm{d}\Omega.
%	\emph{total helicity (or simply helicity) } & \mathcal{H} := \int_{\Omega} \bu \cdot \omega \ \mathrm{d}\Omega.
\end{aligned}
\]
In the absence of the external force\footnote{Or when the external force is conservative.}, i.e., $ \boldsymbol{f}=\boldsymbol{0} $, and if there is no net flux of energy and enstrophy over the domain boundary, two-dimensional incompressible flows dissipate energy at a rate,
\begin{equation}\label{Eq: strong balance a}
	\partial_{t}\mathcal{K} = -2 \Rn^{-1} \mathcal{E},
\end{equation}
and, additionally dissipate enstrophy at a rate,
\begin{equation}\label{Eq: strong balance b}
	\partial_{t}\mathcal{E} = -2 \Rn^{-1} \mathcal{P} .
%	\partial_{t}\mathcal{W} &= \Rn^{-1}\oint_{\partial\Omega}\left( \nabla\times\omega\right) \times\bn\ \mathrm{d}\Gamma + \oint_{\partial\Omega}\left( \omega\times\bu\right) \times\bn\ \mathrm{d}\Gamma + \int_{\Omega}\nabla\times\boldsymbol{f}\ \mathrm{d}\Omega, \label{Eq: strong balance c}
\end{equation}
Furthermore, if the two-dimensional incompressible flow is in the inviscid limit ($ \Rn\to\infty $), namely, the flow is ideal, from \eqref{Eq: strong balance a} and \eqref{Eq: strong balance b}, it conserves energy and enstrophy over time \cite{Palha2017}.
%Note that, in two-dimensions,
%% equipped with a canonical coordinate system
%we can express velocity (or the unit outward norm vector) as a three-dimensional vector whose third component perpendicular to the plane is constantly zero, i.e., $ \bu=\begin{bmatrix}
%	u & v & 0
%\end{bmatrix}^\mathsf{T} $ (or $\bn=\begin{bmatrix}
%n_{x} & n_{y} & 0
%\end{bmatrix}^\mathsf{T} $), and express vorticity as a three-dimensional vector whose in-plane components are both zero, i.e., \( \omega = \begin{bmatrix}0& 0&  \omega\end{bmatrix}^\mathsf{T}\).  Thus relation \eqref{Eq: NSE b}, also written as $ \omega = \nabla\times\bu $ in two-dimensions, implies
And the relation \eqref{Eq: NSE b} implies
\begin{equation}\label{Eq: total vorticity}
	\mathcal{W} := \int_{\Omega}\omega\ \mathrm{d}\Omega = \oint_{\partial\Omega}\bu\times\bn\ \mathrm{d}\Gamma,
\end{equation}
where $ \mathcal{W} $ is called the total vorticity. \eqref{Eq: total vorticity} shows that the total vorticity no matter whether the flow is ideal or not, is a conserved quantity over time provided $ \oint_{\partial\Omega}\bu\times\bn\ \mathrm{d}\Gamma $ is not time dependent.

The first scheme that is mass, energy, enstrophy and vorticity conserving (MEEVC) was proposed in \cite{Palha2017} where two evolution equations for velocity and vorticity are employed. The two evolution equations are staggered in time such that information can be transferred between each other through a midpoint temporal discretization scheme. As a result, both equations are linearized and the unknowns are decoupled to separate time instant sequences, which significantly lowers the computational cost. A drawback of this scheme is that the inclusion of no-slip boundary conditions requires indirect approaches and the suggested approach destroys the vorticity conservation property \cite{DeDiego2019}. An extensive literature study on structure-preserving methods is given in \cite{Palha2017}. For a more recent discussion on structure-preserving methods for incompressible flows, we refer, for example, to \cite{GAWLIK2020109439}.
%For a literature survey on energy-conserving temporal discretization methods for incompressible flows, we refer to \cite{sanderse2013energy}.

In this work, we introduce a mixed high-order finite element discretization of two-dimensional incompressible Navier-Stokes equations that is also MEEVC, can incorporate general boundary condition but avoids the evolution equation for vorticity. The functional setting will be given in Section~\ref{Sec formulation}. In Section~\ref{Sec semi discrete}, properties of the formulation are analyzed, which is followed by the introduction of the temporal discretization in Section~\ref{Sec temporal discrete}. Numerical tests are presented in Section~\ref{Sec numerical tests}. Finally, conclusions are drawn in Section~\ref{Sec conclusions}.

%_____________________________________________________________________________________________________
%>>>>>>>>>>>>>>>>>>>>>>>>>>>>>>>>>>>>>> SECTION <<<<<<<<<<<<<<<<<<<<<<<<<<<<<<<<<<<<<<<<<<<<<<<<<<<<<<<
\section{The mixed weak formulation}\label{Sec formulation}
%Without loss of generality, in this section, we work with boundary conditions
%\[
%\widehat{P}=0,\ \bu\times\bn = \boldsymbol{0}\quad \text{on}\ \partial\Omega\times(0,T],
%\]
%for simplicity. 
\subsection{A brief introduction to function spaces employed}
The space of square integrable functions, is
\[
L^2(\Omega) := \left\lbrace f\ \left|\  \left\langle f,f\right\rangle_{\Omega} \leq +\infty\right.\right\rbrace ,
\]
where $ \left\langle \cdot,\cdot\right\rangle_{\Omega} $ denotes the $ L^2 $-inner product (or simply inner product) over the domain $ \Omega $.
In $\mathbb{R}^2$, we will also use Sobolev spaces
\[
H(\mathrm{curl};\Omega) := \left\lbrace \omega \left| \omega\in  L^2(\Omega),\ \nabla\times\omega \in \left[ L^2(\Omega)\right] ^2 \right.\right\rbrace,
\]
\[
H(\mathrm{div};\Omega) :=\left\lbrace \boldsymbol{u}\left| \boldsymbol{u}\in\left[  L^2(\Omega)\right] ^2,\ \nabla\cdot\boldsymbol{u} \in L^2(\Omega) \right.\right\rbrace.
\]
\[
H^1(\Omega) := \left\lbrace \phi \left| \phi\in  L^2(\Omega),\ \nabla \phi \in \left[ L^2(\Omega)\right] ^2 \right.\right\rbrace,
\]
\[
H(\mathrm{rot};\Omega) :=\left\lbrace \boldsymbol{\sigma}\left| \boldsymbol{\sigma}\in\left[  L^2(\Omega)\right] ^2,\ \nabla\times\boldsymbol{\sigma} \in L^2(\Omega) \right.\right\rbrace.
\]
They form de Rham complexes \cite{10.1007/0-387-38034-5_5} in two dimensions written as
\[
\mathbb{R}\hookrightarrow H(\mathrm{curl};\Omega) \stackrel{\nabla\times}{\longrightarrow}  H(\mathrm{div};\Omega) \stackrel{\nabla\cdot}{\longrightarrow}  L^2(\Omega) \to 0,
\]
\[
\mathbb{R}\hookrightarrow H^1(\Omega)\stackrel{\nabla}{\longrightarrow}  H(\mathrm{rot};\Omega) \stackrel{\nabla\times}{\longrightarrow}  L^2(\Omega) \to 0.
\]
Also see (1) and (2) of \cite{christiansen2018nodal}.

The trace operator, denoted by $\mathcal{T}$, restricts a function to a boundary section, $\Gamma\subseteq\partial\Omega$. The trace operator acting on $\omega\in H(\mathrm{curl};\Omega)$, $\phi \in H^1(\Omega)$, $\bu \in H(\mathrm{div};\Omega)$ and $\boldsymbol{\sigma} \in H(\mathrm{rot};\Omega)$ is, respectively,
\[
\mathcal{T}\omega 
= \left.\omega\right|_{\Gamma} , \quad
\mathcal{T}\phi  = \left.\phi\right|_{\Gamma},
\]
\[
 \mathcal{T}\bu = \left.\bu\cdot\boldsymbol{n}\right|_{\Gamma} ,
\quad \mathcal{T}\boldsymbol{\sigma} =\left.\boldsymbol{\sigma}\times\boldsymbol{n}\right|_{\Gamma} .
\]
%%Recall that $\omega = \begin{bmatrix}
%%	0 & 0 & \omega
%%\end{bmatrix}^{\mathsf{T}}$. 
The trace spaces are 
\[
\mathcal{T}H(\mathrm{curl};\Omega,\Gamma) := \left\lbrace \left.\mathcal{T}\omega\ \right|\ \omega\in H(\mathrm{curl};\Omega)\right\rbrace, 
\]
\[
H^{-1/2}(\Omega,\Gamma) := \left\lbrace \left.\mathcal{T}\bu\ \right|\ \bu\in H(\mathrm{div};\Omega)\right\rbrace, 
\]
\[
H^{1/2}(\Omega,\Gamma) := \left\lbrace \left.\mathcal{T}\phi\ \right|\ \phi\in H^1(\Omega)\right\rbrace, 
\]
\[
\mathcal{T}H(\mathrm{rot};\Omega,\Gamma) := \left\lbrace \left.\mathcal{T}\boldsymbol{\sigma}\ \right|\ \boldsymbol{\sigma}\in H(\mathrm{rot};\Omega)\right\rbrace.
\]
Not that $\mathcal{T}H(\mathrm{curl};\Omega,\Gamma)$ and $\mathcal{T}H(\mathrm{rot};\Omega,\Gamma)$, as well as $H^{-1/2}(\Omega,\Gamma)$ and $H^{1/2}(\Omega,\Gamma)$, are a pair of dual spaces. For a complete introduction on Sobolev spaces, we refer to \cite{Oden2010}.

We use notations
$C(\Omega),\ D(\Omega)$,\ $G(\Omega),\ R(\Omega)$ and $S(\Omega)$
to express  finite dimensional conforming function spaces which are subsets of Sobolev spaces, i.e.,
\[
	 C(\Omega)\subset  H(\mathrm{curl};\Omega),\quad D(\Omega)\subset H(\mathrm{div};\Omega),\quad S(\Omega)\subset L^2(\Omega) ,
\]
\[
G(\Omega)\subset  H^1(\Omega),\quad R(\Omega)\subset H(\mathrm{rot};\Omega),
\]
and form discrete de Rham complexes in two-dimensional space,
\begin{equation}\label{Eq 2d de Rham}
	\mathbb{R}\hookrightarrow C(\Omega)\stackrel{\nabla\times}{\longrightarrow}\ D(\Omega)\stackrel{\nabla\cdot}{\longrightarrow}\ S(\Omega)\to 0.
\end{equation}
\begin{equation*}
	\mathbb{R}\hookrightarrow G(\Omega)\stackrel{\nabla}{\longrightarrow}\ R(\Omega)\stackrel{\nabla\times}{\longrightarrow}\ S(\Omega)\to 0.
\end{equation*}
The finite dimensional spaces $ C(\Omega) $ and $ D(\Omega) $ possess sufficient regularity such that
\begin{equation}\label{Eq: regularity}
\omega_{h}\times\boldsymbol{u}_{h} \in \left[ L^2(\Omega)\right] ^{2},\quad  \forall \left( \omega_{h}, \boldsymbol{u}_{h}\right) \in C(\Omega)\times D(\Omega).
\end{equation} 
Trace spaces of finite dimensional spaces $C(\Omega,\Gamma) $ and $d(\Omega,\Gamma) $ on boundary section $\Gamma$ are denoted by
\[
\mathcal{T}C(\Omega,\Gamma) := \left\lbrace \left.\mathcal{T}\omega_{h}\ \right|\ \omega_{h}\in C(\Omega)\right\rbrace, 
\]
\[
\mathcal{T}D(\Omega,\Gamma) := \left\lbrace \left.\mathcal{T}\boldsymbol{u}_{h}\ \right|\ \boldsymbol{u}_{h}\in D(\Omega)\right\rbrace. 
\]
%or a part of a discrete de Rham complex in $ \mathbb{R}^3 $, i.e.,
%\begin{equation}\label{Eq 3d de Rham}
%	\cdots\stackrel{\nabla}{\longrightarrow} \ C(\Omega)\stackrel{\nabla\times}{\longrightarrow}\ D(\Omega)\stackrel{\nabla\cdot}{\longrightarrow}\ S(\Omega)\to 0.
%\end{equation}
%a subspace of $ H(\mathrm{div};\Omega) $,
%\[
%H_{\widehat{u}}(\mathrm{div};\Omega):=\left\{\bu\left|\bu\in H(\mathrm{div};\Omega),\  T\bu = \widehat{u} \ \text{on}\ \Gamma_{\widehat{u}}\right.\right\},
%\]
And we will also use following subspaces,
\[C_{0}(\Omega,\Gamma) :=\left\lbrace  \omega_{h}\left|\ \omega_{h}\in C(\Omega), \ \mathcal{T} \omega_{h} = 0 \in \mathcal{T}C(\Omega,\Gamma)\right. \right\rbrace ,\]
\[D_{0} (\Omega, \Gamma) :=\left\lbrace  \boldsymbol{u}_{h}\left|\ \boldsymbol{u}_{h}\in D(\Omega), \ \mathcal{N}\boldsymbol{u}_{h} = 0 \in \mathcal{N}D(\Omega,\Gamma)\right. \right\rbrace .\]
%and their trace spaces,
%\[
%\left[ T D_{\chi} (\Omega) \right] _{00} := \left\lbrace T\boldsymbol{b}_{h}\left| \boldsymbol{b}_{h} \in D_{\chi} (\Omega),\ T\boldsymbol{b}_{h} = 0\ \text{on}\ \Gamma_{\widehat{P}} \right.  \right\rbrace ,
%\]
%\[
%\left[ T^{\parallel} C_{\chi} (\Omega) \right] _{00} := \left\lbrace T^{\parallel}\boldsymbol{a}_{h}\left| \boldsymbol{a}_{h} \in C_{\chi} (\Omega),\  T^{\parallel}\boldsymbol{a}_{h} = \boldsymbol{0}\ \text{on}\ \Gamma_{u^{\parallel}} \right.  \right\rbrace .
%\]

\subsection{The formulation}
We introduce a trilinear form
\[
a\left( \rho_{h}, \boldsymbol{\vartheta}_{h}, \boldsymbol{e}_{h} \right)
%:=\left\langle \rho_{h}\times\boldsymbol{\vartheta}_{h},\boldsymbol{\varepsilon}_{h} \right\rangle_{\Omega}  
:= \int_{\Omega} \left(\rho_{h} \times \boldsymbol{\vartheta}_{h}\right) \cdot \boldsymbol{e}_{h} \, \mathrm{d}\Omega,
\]
for $\left( \rho_{h}, \boldsymbol{\vartheta}_{h}, \boldsymbol{e}_{h} \right)\in C(\Omega)\times D(\Omega)\times D(\Omega) $. Because $\rho_{h}\times\boldsymbol{\vartheta}_{h}$ is pointwise perpendicular to $\boldsymbol{\vartheta}_{h}$, we know that 
\begin{equation}\label{Eq: perp}
	a\left( \rho_{h}, \boldsymbol{\vartheta}_{h}, \boldsymbol{\vartheta}_{h} \right)=0.
\end{equation}

A spatially discrete weak mixed formulation of \eqref{Eq: NSE} is written as following: Given
$ \boldsymbol{f} \in \left[ L^{2}(\Omega) \right] ^{2} $
and natural boundary conditions,
$
\widehat{P}\in H^{1/2}(\Omega, \Gamma_{\widehat{P}})$ and $
\widehat{u}_{\parallel}\in \mathcal{T}H(\mathrm{rot};\Omega,\Gamma_{\parallel})
$,
%$ \bu_{h}^{0}\in H_{\widehat{u}}(\mathrm{div};\Omega)$ 
seek $ \left(\bu_{h},\omega_{h},P_{h}\right) \in D(\Omega)\times C(\Omega)\times S(\Omega)$, such that, $ \forall \left( \bv_{h},\ \xi_{h},\ q_{h}\right) \in D_{0}(\Omega,\Gamma_{\perp})\times  C_{0}(\Omega,\Gamma_{\widehat{\omega}})\times S(\Omega)$,   
\begin{subequations}\label{Eq: form 1}
	\begin{align}
		\left\langle \partial_{t}\bu_{h}, \bv_{h}\right\rangle_{\Omega}  + 
		a\left(  \omega_{h},\bu_{h}, \bv_{h}\right)  +\Rn^{-1}\left\langle\nabla\times\omega_{h}, \bv_{h}\right\rangle_{\Omega} - \left\langle  P_{h} , \nabla\cdot\bv_{h}\right\rangle_{\Omega} &= \left\langle \bF, \bv_{h}\right\rangle_{\Omega} - \left\langle\left.\widehat{P}\right|\mathcal{T}\bv_{h} \right\rangle_{\Gamma_{\widehat{P}}} ,\label{Eq: a}\\
		\left\langle  \bu_{h}, \nabla\times\xi_{h} \right\rangle_{\Omega} -\left\langle \omega_{h} , \xi_{h}\right\rangle_{\Omega}&=\left\langle \left.\widehat{u}_{\parallel}\right|\mathcal{T}\xi_{h} \right\rangle_{\Gamma_{\parallel}} ,\label{Eq: b}\\
		\left\langle \nabla\cdot\bu_{h},q_{h}\right\rangle_{\Omega} &= 0,\label{Eq: c}
	\end{align}
\end{subequations}
subject to essential boundary conditions, $
\mathcal{T}\boldsymbol{u}_{h} = \widehat{u}_{\perp} \in \mathcal{T} D(\Omega,\Gamma_{\perp})$ and $
\mathcal{T} \omega_{h} = \widehat{\omega} \in \mathcal{T} C(\Omega, \Gamma_{\widehat{\omega}})$, and initial conditions $\left(\bu_{h}^{0},\omega_{h}^{0}\right)\in D(\Omega)\times C(\Omega)$. Note that we have used the notation $\left<\left.\cdot\right|\cdot\right>_{\Gamma}$ to indicate that it is a duality pairing between elements from a pair of dual spaces.

One can show that $ \left(\bu_{h},\omega_{h},P_h\right) $ in \eqref{Eq: form 1} weakly solves the incompressible Navier-Stokes equations \eqref{Eq: NSE}, see the work of Boffi, Brezzi and Fortin \cite{boffi2013mixed}. A similar setup is proposed in the fluid part of a structure-preserving formulation for magnetohydrodynamics (MHD) problems, cf. (48) - (53) of \cite{GAWLIK2022110847}.

\section{Dissipation and conservation properties} \label{Sec semi discrete}
In this section, we study dissipation and conservation properties of the formulation \eqref{Eq: form 1}. The equivalence between \eqref{Eq: form 1} and the formulation used in the original MEEVC scheme, see (19) of \cite{Palha2017}, will also be shown. Thus, we will prove that the formulation \eqref{Eq: form 1} is also MEEVC. To this end, the analysis conducted here is under conditions that (i) the domain is periodic ($ \partial\Omega = \emptyset $) and (ii) there is no external force as in \cite{Palha2017}.

\subsection{Mass conservation}
Pointwise mass conservation is obviously satisfied; $ \bu_{h} $ is selected to be in $D(\Omega)\subset H(\mathrm{div};\Omega) $ and the relation \eqref{Eq: c} strongly enforces $ \nabla\cdot\bu_{h} = 0 $ everywhere in $ \Omega $. This is a consequence of the fact that $\nabla\cdot$ maps $D(\Omega)$ into $S(\Omega)$, see \eqref{Eq 2d de Rham}.

\subsection{Energy dissipation and conservation}
For the energy balance, if we replace $ \bv_{h} $ in \eqref{Eq: a} by $ \bu_{h}\in D(\Omega) $, we will obtain
\begin{equation}\label{Eq: ec0}
	\left\langle \partial_{t}\bu_{h}, \bu_{h}\right\rangle_{\Omega} 
	+ a\left(  \omega_{h},\bu_{h}, \bu_{h}\right) 
	+\Rn^{-1}\left\langle\nabla\times\omega_{h}, \bu_{h}\right\rangle_{\Omega}
	- \left\langle  P_{h} , \nabla\cdot\bu_{h}\right\rangle_{\Omega} 
	= 0
	%\left\langle \bF, \bu_{h}\right\rangle_{\Omega} 
	%- \left\langle\widehat{P},T\bu_{h} \right\rangle_{\partial\Omega}
	.
\end{equation}
The second and fourth terms vanish because of \eqref{Eq: perp} and the pointwise mass conservation, i.e., $ \nabla\cdot\bu_{h} = 0 $, respectively. This leads to
\begin{equation*}
	\left\langle \partial_{t}\bu_{h}, \bu_{h}\right\rangle_{\Omega} 
	+ \Rn^{-1}\left\langle\nabla\times\omega_{h}, \bu_{h}\right\rangle_{\Omega}
	= 0
	%\left\langle \bF, \bu_{h}\right\rangle_{\Omega} 
	%- \left\langle\widehat{P},T\bu_{h} \right\rangle_{\partial\Omega}
	.
\end{equation*}
And from \eqref{Eq: b}, we know
\begin{equation*}
	\left\langle  \bu_{h}, \nabla\times{\omega}_{h} \right\rangle_{\Omega}
	= \left\langle \omega_{h} , {\omega}_{h}\right\rangle_{\Omega} 
	%+\left\langle u^{\parallel},T^{\parallel}\boldsymbol{\omega}_{h} \right\rangle_{\partial\Omega}
	,
\end{equation*}
because $ \omega_{h}\in C(\Omega) $. Combining these two relations 
%and rearranging the terms will 
gives a (semi-)discrete energy balance,
\begin{equation}\label{Eq: ec1}
		\partial_{t}\mathcal{K}_{h} = \left\langle \partial_{t}\bu_{h}, \bu_{h}\right\rangle_{\Omega} 
	= - \Rn^{-1}\left\langle \omega_{h} , {\omega}_{h}\right\rangle_{\Omega}  
	=- 2\Rn^{-1}\mathcal{E}_{h} 
	%-\Rn^{-1}\left\langle u^{\parallel},T^{\parallel}\boldsymbol{\omega}_{h} \right\rangle_{\partial\Omega} - \left\langle\widehat{P},T\bu_{h} \right\rangle_{\partial\Omega} + \left\langle \bF, \bu_{h}\right\rangle_{\Omega}
	,
\end{equation}
where $ \mathcal{K}_{h} = \dfrac{1}{2}\left\langle\bu_{h},\bu_{h}\right\rangle_{\Omega}$ and $ \mathcal{E}_{h} = \dfrac{1}{2}\left\langle\omega_{h},\omega_{h}\right\rangle_{\Omega} $ are the discrete (total kinetic) energy and (total) enstrophy, respectively. It is consistent with \eqref{Eq: strong balance a}, the energy balance of the strong form. Thus, \eqref{Eq: ec1} clearly implies discrete energy conservation in the inviscid limit ($\Rn\to\infty$).

\subsection{Enstrophy dissipation and conservation}
If we take the time derivative of \eqref{Eq: b}, we obtain
\begin{equation}\label{Eq: helicity conservation relation 1}
	\left\langle\partial_t\bu_{h}, \nabla\times \bx_{h}\right\rangle_{\Omega} = \left\langle \partial_t\omega_{h},\bx_{h} \right\rangle_{\Omega}, \quad \forall\bx_{h}\in C(\Omega)\,.
\end{equation}
And, from \eqref{Eq: a}, we know that, $\forall\ \bx_{h} \in C(\Omega) $, \eqref{Eq: a} must hold for $ \nabla\times\bx_{h}\in D(\Omega) $, namely,
\begin{equation}\label{Eq: enstrophy 0}
	\left\langle\partial_t\bu_{h}, \nabla\times\bx_{h}\right\rangle_{\Omega}  + 
	a\left(  \omega_{h},\bu_{h}, \nabla\times\bx_{h}\right)  +\Rn^{-1}\left\langle\nabla\times\omega_{h}, \nabla\times\bx_{h}\right\rangle_{\Omega} - \left\langle  P_{h} , \nabla\cdot\nabla\times\bx_{h}\right\rangle_{\Omega} =  0
	%\left\langle \bF, \nabla\times\xi_{h}\right\rangle_{\Omega} 
	,
\end{equation} 
where the term $ \left\langle  P_{h} , \nabla\cdot\nabla\times\bx_{h}\right\rangle_{\Omega} $ 
%and $ \left\langle\widehat{P},T\left( \nabla\times\xi_{h} \right) \right\rangle_{\partial\Omega} $ 
vanishes because $ \nabla\cdot\nabla\times(\cdot)\equiv 0 $.
%, and the later vanishes because $ T\left( \nabla\times\xi_{h} \right)= \left( \nabla\times\xi_{h} \right)\cdot\bn = 0 $ as $ \left( \nabla\times\xi_{h}\right) \perp \bn$. 
If we further insert \eqref{Eq: helicity conservation relation 1} into \eqref{Eq: enstrophy 0}, we obtain
\begin{equation}\label{Eq: en-2}
	\left\langle \partial_t\omega_{h},\bx_{h} \right\rangle_{\Omega} + a\left(  \omega_{h},\bu_{h}, \nabla\times\bx_{h}\right)   +\Rn^{-1}\left\langle\nabla\times\omega_{h}, \nabla\times\bx_{h}\right\rangle_{\Omega} =  0,\quad \forall\bx_{h}\in C(\Omega)\,.
\end{equation}
We can replace $\bx_{h}$ in \eqref{Eq: en-2} by $\omega_{h}\in C(\Omega)$ and get
 \begin{equation}\label{Eq: en-3}
 	\left\langle \partial_t\omega_{h},\omega_{h} \right\rangle_{\Omega} +a\left(  \omega_{h},\bu_{h}, \nabla\times\omega_{h}\right)   +\Rn^{-1}\left\langle\nabla\times\omega_{h}, \nabla\times\omega_{h}\right\rangle_{\Omega} =  0 \,.
 \end{equation}
As $\bu_{h}\in D(\Omega)$ and $\nabla\cdot \boldsymbol{u}_{h} = 0 $ is satisfied pointwise, we can find a stream function $\psi_{h}\in C(\Omega)$ (on the simply connected, contractible domain) such that $\bu_{h} = \nabla\times\psi_{h}$. Recall the following vector calculus identity, 
\[\omega_{h} \times \nabla\times\psi_{h} = \nabla\left( \omega_{h}\psi_{h}\right) - \psi_{h}\times \nabla\times \omega_{h}.\]
Thus we know
\begin{equation*}
	\begin{aligned}
		a\left(  \omega_{h},\bu_{h}, \nabla\times\omega_{h}\right)&= 
		\int_{\Omega}\nabla\left( \omega_{h}\psi_{h}\right)\cdot \left( \nabla\times\omega_{h}\right) \mathrm{d}\Omega
		- 
		a\left( \psi_{h},  \nabla\times\omega_{h}, \nabla\times\omega_{h}\right) \\
		& = \int_{\Omega}\omega_{h}\psi_{h} \left( \nabla\cdot\nabla\times\omega_{h}\right) \mathrm{d}\Omega
		- 
		a\left( \psi_{h},  \nabla\times\omega_{h}, \nabla\times\omega_{h}\right),
	\end{aligned}
\end{equation*}
where we have performed integration by parts with respect to the gradient operator for the first term of the second equality and use the periodic boundary condition. Obviously, these terms vanish because of property $\nabla\cdot\nabla\times\left( \cdot\right) \equiv  0$ and \eqref{Eq: perp}. Therefore, we know that 
\begin{equation}\label{Eq: key eq}
	a\left(  \omega_{h},\bu_{h}, \nabla\times\omega_{h}\right) = 0, 
\end{equation}
and \eqref{Eq: en-3} leads to the following (semi-)discrete enstrophy balance,
\begin{equation} \label{Eq: en-4}
	\left\langle \partial_t\omega_{h},\omega_{h} \right\rangle_{\Omega}  =-\Rn^{-1}\left\langle\nabla\times\omega_{h}, \nabla\times\omega_{h}\right\rangle_{\Omega}
	=-2\Rn^{-1}\mathcal{P}_{h}\,,
\end{equation}
where $ \mathcal{P}_{h} := \dfrac{1}{2}\left\langle \nabla\times \omega_{h},\nabla\times\omega_{h}\right\rangle_{\Omega} $ is the discrete (total) palinstrophy.
\eqref{Eq: en-4} correctly reflects the enstrophy balance of the strong form, see \eqref{Eq: strong balance b}. And, in the inviscid limit ($\Rn\to\infty$), \eqref{Eq: en-4} leads to enstrophy conservation. 

%and the conclusion of enstrophy conservation (and dissipation), i.e., \eqref{Eq: enstrophy balance}, follows if we select $\xi_{h}$ to be $\omega_{h}$ in \eqref{Eq: enstrophy 0} . In Section~\ref{Sub: ec}, we will also give numerical evidences showing that \eqref{Eq: eq-new-1} is valid for no matter exact or inexact numerical quadrature. 

\subsection{Vorticity conservation}

For conservation of (total) vorticity, if we select $ \xi_{h} = 1 $ in \eqref{Eq: enstrophy 0}, it is straightforward to find that
\[
\partial_{t}\mathcal{W}_{h} = \left\langle \partial_{t} \omega_{h}, 1\right\rangle_{\Omega} = 0,
\]
which implies that vorticity is conserved over time. Moreover, by selecting $ \bx_{h}=1 $ in \eqref{Eq: b}, we know that in periodic domains
\[
\mathcal{W}_{h} \equiv 0 ,
\]
which is consistent with \eqref{Eq: total vorticity} of the strong form.

\subsection{Equivalence to the original MEEVC formulation}
If we apply  integration by parts to the second term of \eqref{Eq: en-2}, we obtain
\begin{equation}\label{Eq: enstrophy 1}
	\left\langle \partial_{t}\omega_{h} , \bx_{h}\right\rangle_{\Omega} 
	+ \left\langle\nabla\times\left( \omega_{h}\times \bu_{h}\right) , \bx_{h}\right\rangle_{h}
	+\Rn^{-1}\left\langle\nabla\times\omega_{h}, \nabla\times\bx_{h}\right\rangle_{\Omega} = 0.
	% \left\langle \bF, \nabla\times\xi_{h}\right\rangle_{\Omega} - \left\langle\widehat{P},T\left( \nabla\times\xi_{h} \right) \right\rangle_{\partial\Omega},\quad  \forall\xi_{h}\in C(\Omega),
\end{equation}
Recall that the following identity
\begin{equation}\label{Eq:magic-equation-0}
	\nabla\times\left( \omega_{h}\times \bu_{h}\right) = \dfrac{1}{2}\left( \bu_{h}\cdot\nabla\right) \omega_{h} + \dfrac{1}{2}\nabla\cdot\left( \omega_{h}\bu_{h}\right) ,
\end{equation}
is valid in two-dimensions.
Using this identity, the following substitution can be employed in the weak form of the vorticity evolution equation \eqref{Eq: enstrophy 1},
\begin{equation}\label{Eq:magic-equation}
	\begin{aligned}
	\left\langle \nabla\times\left( \omega_{h}\times \bu_{h}\right),\xi_{h}\right\rangle_{\Omega} &= \dfrac{1}{2}\left\langle \left( \bu_{h}\cdot\nabla\right) \omega_{h},\xi_{h}\right\rangle_{\Omega} + 
	\dfrac{1}{2}\left\langle \nabla\cdot\left( \omega_{h}\bu_{h}\right),\xi_{h}\right\rangle_{\Omega}
	\\&= 
	-\dfrac{1}{2}\left\langle \omega_{h}, \nabla\cdot\left(\xi_{h}\bu_{h} \right) \right\rangle_{\Omega} + 
	\dfrac{1}{2}\left\langle \nabla\cdot\left( \omega_{h}\bu_{h}\right),\xi_{h}\right\rangle_{\Omega},
\end{aligned}
\end{equation}
where the following integration by parts was used,
\[
\left\langle \left( \bu_{h}\cdot\nabla\right) \omega_{h},\xi_{h}\right\rangle_{\Omega} = -\left\langle \omega_{h}, \nabla\cdot\left(\xi_{h}\bu_{h} \right) \right\rangle_{\Omega}.
\]
Thus, $ \forall \xi_{h}\in C(\Omega) $, \eqref{Eq: enstrophy 1} can be written as
\begin{equation}\label{Eq: MEEVC}
	\left\langle \partial_{t} \omega_{h}, \xi_{h}\right\rangle_{\Omega}   -\dfrac{1}{2}\left\langle \omega_{h}, \nabla\cdot\left(\xi_{h}\bu_{h}\right) \right\rangle_{\Omega}  + \dfrac{1}{2}\left\langle \nabla\cdot\left( \omega_{h}\bu_{h}\right),\xi_{h} \right\rangle_{\Omega} + \Rn^{-1}\left\langle\nabla\times\omega_{h}, \nabla\times\xi_{h}\right\rangle_{\Omega} = 0,
\end{equation}
which is the weak evolution equation for vorticity (as a replacement of \eqref{Eq: b}) in the original MEEVC scheme, see (19) of \cite{Palha2017}. 
%This concludes the equivalence between formulation \eqref{Eq: form 1} and the formulation used in the original MEEVC scheme (except that in (19) of \cite{Palha2017} particular finite dimensional spaces have been selected). 
%Thus, the proof of enstrophy and total vorticity conservation can be done using a similar method as in \cite{Palha2017}. Just for completeness, we repeat it briefly below. \Yi{It is worthy of mentioning that, to show this equivalence, we have used some steps, see \eqref{Eq: enstrophy 1} and \eqref{Eq:magic-equation}, which are valid only under exact numerical intergration.  But it does not mean that the present method works only for exact numerical intergration. In Section~\ref{subsec SLR}, we will see numerical evidences that the present method is still MEEVC in an extremely distorted mesh where exact numerical intergration is not employed. Analytically, this remains an open question to the authors.}
If we select $ \xi_{h} $ to be $ \omega_{h} $ in \eqref{Eq: MEEVC}, the second and third terms cancel, we again get the same enstrophy balance, i.e.,\eqref{Eq: en-4}. 

In the original MEEVC work, the reason behind using the identity \eqref{Eq:magic-equation-0} is to replace the weak nonlinear advection term by \eqref{Eq:magic-equation}. By doing so, it is possible to construct another trilinear form for the nonlinear advection term in the weak vorticity evolution equation which is skew-symmetric with respect to entries $\omega_{h}$ and $\xi_{h}$, i.e.,
\[
b(\omega_{h}, \bu_{h},\xi_{h}) = - b(\xi_{h}, \bu_{h},\omega_{h}).
\]
See second and third terms in \eqref{Eq: MEEVC}. This then implies that 
\[
b(\omega_{h},\bu_{h},\omega_{h}) = 0,
\]
which, for the original MEEVC work, is a key requirement to obtain enstrophy conservation even when the numerical quadrature is inexact. However, this prevents the direct incorporation of boundary conditions for the tangential component of velocity because the port, i.e., the boundary integral term in \eqref{Eq: b} which is used to impose them weakly, is missing. In this present work, we surprisingly find that, to setup a MEEVC scheme, we can bypass the construction of the skew-symmetric advection term as in the original MEEVC scheme and, thus, it is not necessary to introduce a second evolution equation for vorticity. This simplifies the formulation and also enables the direct application of no-slip boundary conditions.

%What is surprising in the present work is that the original term is skew-symmetric at the discrete level, i.e.,
%\[
%\left\langle \omega_{h}\times\bu_{h} ,\nabla\times\xi_{h} \right\rangle _{\Omega}
%=-
%\left\langle \xi_{h}\times\bu_{h} ,\nabla\times\omega_{h} \right\rangle _{\Omega}.
%\]
%With exact quadrature this identity is straightforward, but it is not when quadrature is not exact (which usually is the case on curved geometries). We have numerically tested this properly and have verified that enstrophy conservation is indeed satisfied when the quadrature is inexact, cf. Section~\ref{subsec SLR}. We intend to analytically prove this in the future and understand if this is dependent on the basis functions used (e.g., the mimetic spectral elements \cite{Kreeft2011,zhang2022phd} used in this work versus classic mixed finite elements \cite{boffi2013mixed}.

\section{Temporal discretization} \label{Sec temporal discrete}
For the temporal discretization, the classic implicit midpoint method \cite{hairer2006geometric} is used. The fully discrete version of \eqref{Eq: form 1} is written as: Given
$ \boldsymbol{f} \in \left[ L^{2}(\Omega) \right] ^{2} $
and natural boundary conditions,
$
\widehat{P}\in H^{1/2}(\Omega, \Gamma_{\widehat{P}})$ and $
\widehat{u}_{\parallel}\in \mathcal{T}H(\mathrm{rot};\Omega,\Gamma_{\parallel})
$,
%$ \bu_{h}^{0}\in H_{\widehat{u}}(\mathrm{div};\Omega)$ 
for $ k\in\left\{1,2,3,\cdots\right\} $,  seek $ \left(\bu^{k}_{h},\omega^{k}_{h},P^{k-\frac{1}{2}}_{h}\right) \in D(\Omega)\times C(\Omega)\times S(\Omega)$, such that, $ \forall \left( \bv_{h},\ \xi_{h},\ q_{h}\right) \in D_{0}(\Omega,\Gamma_{\perp})\times  C_{0}(\Omega,\Gamma_{\widehat{\omega}})\times S(\Omega)$, 
\begin{subequations}\label{Eq fully discrete}
	\begin{align}
		&\left\langle \frac{\bu^{k}_{h}-\bu_{h}^{k-1}}{\Delta t}, \bv_{h}\right\rangle_{\Omega}  
		+ 
		a\left( \dfrac{\omega^{k-1}_{h}+\omega^{k}_{h}}{2},\dfrac{\bu^{k-1}_{h}+\bu^{k}_{h}}{2}, \bv_{h}\right)
		 +\Rn^{-1}\left\langle\nabla\times\dfrac{\omega^{k-1}_{h}+\omega^{k}_{h}}{2}, \bv_{h}\right\rangle_{\Omega} \\ 
		&\nonumber\hspace{6cm}- \left\langle  P^{k-\frac{1}{2}}_{h} , \nabla\cdot\bv_{h}\right\rangle_{\Omega} =
		\left\langle \bF^{k-\frac{1}{2}}, \bv_{h}\right\rangle_{\Omega} - \left\langle\left.\widehat{P}^{k-\frac{1}{2}}\right|\mathcal{T}\bv_{h} \right\rangle_{\Gamma_{\widehat{P}}} ,\\
		&\left\langle  \bu^{k}_{h}, \nabla\times\xi_{h} \right\rangle_{\Omega} -\left\langle \omega^{k}_{h} , \xi_{h}\right\rangle_{\Omega}=\left\langle\left. \widehat{u}_{\parallel}^{k}\right|\mathcal{T}\xi_{h} \right\rangle_{\Gamma_{\parallel}} ,\label{Eq fully discrete b}\\
		&\left\langle \nabla\cdot\bu_{h}^{k},q_{h}\right\rangle_{\Omega} = 0, \label{Eq fully discrete c}
	\end{align}
\end{subequations}
where $ \Delta t = t_{k} - t_{k-1} > 0 $, 
%$ \omega_{h}^{0} = \nabla\times\bu_{h}^{0} $, 
$ \bu_{h}^{k} = \bu_{h}(\boldsymbol{x}, t_{k}) $ (see \eqref{Eq initial conditon}), subject to essential boundary conditions, 
$
\mathcal{T}\boldsymbol{u}^{k}_{h} = \widehat{u}_{\perp} \in \mathcal{T} D(\Omega,\Gamma_{\perp})$ and $
\mathcal{T} \omega^{k}_{h} = \widehat{\omega} \in \mathcal{T} C(\Omega, \Gamma_{\widehat{\omega}})$, and initial conditions $\left(\bu_{h}^{0},\omega_{h}^{0}\right)\in D(\Omega)\times C(\Omega)$.

At the fully discrete level, if we repeat the analysis in Section~\ref{Sec semi discrete} now for the fully discrete formulation \eqref{Eq fully discrete}, we can find that pointwise conservation of mass is satisfied at each time instant, see \eqref{Eq fully discrete c}, i.e.,
\[
\nabla\cdot\bu_{h}^{k} = 0
\]
everywhere in $ \Omega $. And we can also obtain dissipation rates,
\begin{equation}\label{Eq discrete energy dissipation}
	\begin{aligned}
		\dfrac{\mathcal{K}^{k}_{h}-\mathcal{K}^{k-1}_{h}}{\Delta t} 
		= \left\langle \frac{\bu^{k}_{h}-\bu_{h}^{k-1}}{\Delta t}, \frac{\bu^{k-1}_{h}+\bu_{h}^{k}}{2}\right\rangle_{\Omega} 
%		&= - \Rn^{-1}\left\langle \frac{\omega^{k-1}_{h}+\omega_{h}^{k}}{2}, \frac{\omega^{k-1}_{h}+\omega_{h}^{k}}{2}\right\rangle_{\Omega} \\
%		&\hspace{0cm}=
=
		- \Rn^{-1}\left\langle \omega_{h}^{k-\frac{1}{2}}, \omega_{h}^{k-\frac{1}{2}}\right\rangle_{\Omega}
		= -2\Rn^{-1}\mathcal{E}_{h}^{k-\frac{1}{2}},
	\end{aligned}
\end{equation}
\begin{equation}\label{Eq discrete enstrophy dissipation}
	\begin{aligned}
		\dfrac{\mathcal{E}^{k}_{h}-\mathcal{E}^{k-1}_{h}}{\Delta t} 
		= \left\langle \frac{\omega^{k}_{h}-\omega_{h}^{k-1}}{\Delta t}, \frac{\omega^{k-1}_{h}+\omega_{h}^{k}}{2}\right\rangle_{\Omega} 
%		&= - \Rn^{-1}\left\langle \nabla\times\frac{\omega^{k-1}_{h}+\omega_{h}^{k}}{2}, \nabla\times\frac{\omega^{k-1}_{h}+\omega_{h}^{k}}{2}\right\rangle_{\Omega} \\
%		&\hspace{0cm}= 
=
		- \Rn^{-1}\left\langle \nabla\times\omega_{h}^{k-\frac{1}{2}}, \nabla\times\omega_{h}^{k-\frac{1}{2}}\right\rangle_{\Omega}
		=-2\Rn^{-1}\mathcal{P}_{h}^{k-\frac{1}{2}},
	\end{aligned}
\end{equation}
and
\begin{equation*} %\label{Eq discrete vorticity dissipation}
	\dfrac{\mathcal{W}^{k}_{h}-\mathcal{W}^{k-1}_{h}}{\Delta t} = 
	\left\langle \frac{\omega^{k}_{h}-\omega_{h}^{k-1}}{\Delta t}, 1\right\rangle_{\Omega} = 0,
\end{equation*}
where 
$\omega_{h}^{k-\frac{1}{2}}:=\frac{\omega^{k-1}_{h}+\omega_{h}^{k}}{2}$.
%where $ \omega_{h}^{k+\frac{1}{2}} = \dfrac{\omega^{k+1}_{h}+\omega_{h}^{k}}{2} $ and $ \omega_{h}^{k+\frac{1}{2}} = \dfrac{\omega^{k+1}_{h}+\omega_{h}^{k}}{2} $.
%Note that, \eqref{Eq discrete enstrophy dissipation} and \eqref{Eq discrete vorticity dissipation} are valid only in two-dimensions.
%, and, because we compute $ \omega_{h}^{0} $ through $ \omega_{h}^{0} = \nabla\times\bu_{h}^{0} $ instead of \eqref{Eq fully discrete b}, \eqref{Eq discrete enstrophy dissipation} is valid for $ k>0 $. 
This shows that at the fully discrete level the proposed scheme is also MEEVC in the inviscid limit, $\Rn\to\infty$.

\section{Numerical tests}\label{Sec numerical tests}
Four tests are conducted in this work. The accuracy of the method is investigated with an analytical solution in Section~\ref{subsec TGV}. Conservation and dissipation properties are tested in Section~\ref{subsec SLR}. The original MEEVC scheme experiences difficulties of handling no-slip boundary conditions \cite{DeDiego2019} while imposing different boundary conditions including the no-slip ones is straightforward for the method studied in this work. This is demonstrated in Section~\ref{subsec NDC}.  Numerical evidences of \eqref{Eq: key eq} is given by the test in Section~\ref{Sub: ec}.

\begin{figure}[h!]
	\centering{
		\subfloat{
			\begin{minipage}[b]{0.5\textwidth}
				\centering
				\includegraphics[width=0.55\linewidth]{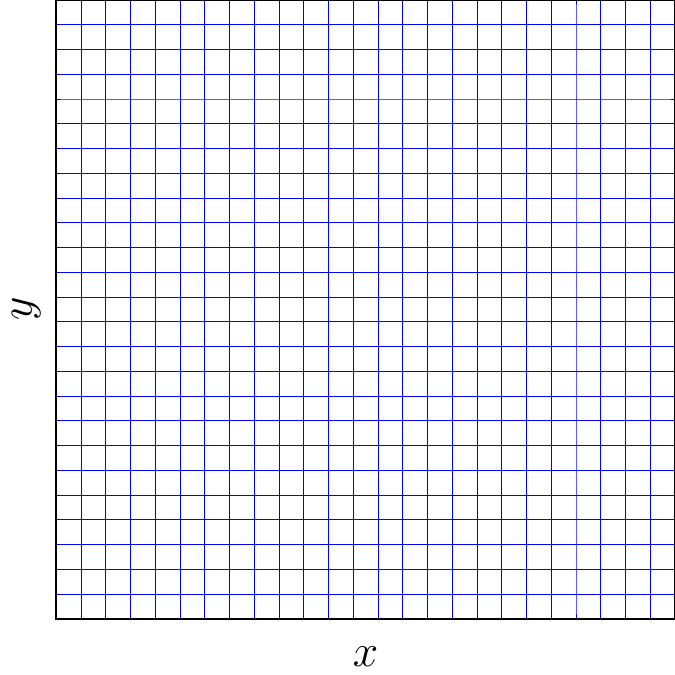}
			\end{minipage}
		}
		\subfloat{
			\begin{minipage}[b]{0.5\textwidth}
				\centering
				\includegraphics[width=0.55\linewidth]{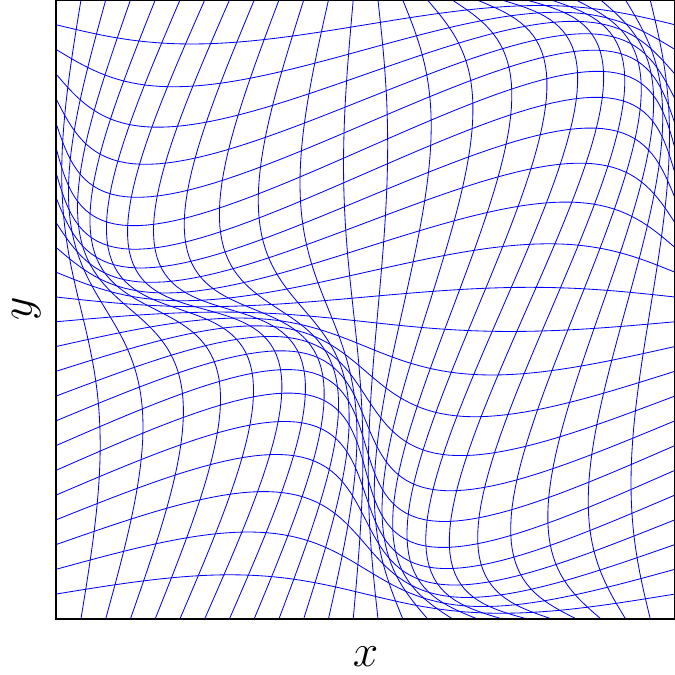}
			\end{minipage}
		}
		\caption{Illustrations of meshes implied by \eqref{eq:crazymesh} for $K=25$, deformation factor $c=0$ (left) and $c=0.25$ (right).}
		\label{fig:crazymesh}
	}
\end{figure}

We use the mimetic polynomial spaces, which satisfy the discrete de Rham complex \eqref{Eq 2d de Rham} and the regularity \eqref{Eq: regularity}, as the finite dimensional spaces under the framework of the mimetic spectral element method, see, for example, \cite{Kreeft2011} or \cite[Chapter~2]{zhang2022phd}. The degree of the polynomial spaces is denoted by $ N $. The Newton-Raphson method is employed for solving the nonlinear systems. Both orthogonal and curvilinear meshes will be used. Suppose a reference domain is $\mathring{\Omega}:=\left(r,s\right)\in[0,1]^2$. A uniform orthogonal mesh of $K\times K$ square elements is generated in the reference domain. This mesh is then distorted with a mapping, $\Phi:\left(r,s\right)\to \left(x,y\right)$, expressed as
\begin{equation}\label{eq:crazymesh}
\left\lbrace
\begin{aligned}
	&x = \alpha \left( r + \frac{1}{2}c\sin(2\pi r)\sin(2\pi s) \right) \\
	&y = \alpha \left( s + \frac{1}{2}c\sin(2\pi r)\sin(2\pi s)\right)
\end{aligned},
\right.
\end{equation}
where $\alpha > 0$ and $0\leq c\leq 0.3$. It gives a mesh in $\Omega=(x,y)\in[0,\alpha]^2$, and the factor $c$ is a deformation factor. When $c=0$ the mesh is orthogonal and uniform, and when $c>0$ the mesh is curvilinear. See Fig.~\ref{fig:crazymesh} for illustrations of this mesh. And see \cite{knupp2020fundamentals} for an introduction on mesh deformation.

Implementations of the present work are done in \emph{Python}. 

\subsection{Accuracy test: Taylor–Green vortex} \label{subsec TGV}
\begin{figure}[!htb]
	\centering
	\begin{minipage}[c]{0.125\textwidth}
		\centering
		\includegraphics[width=1\linewidth]{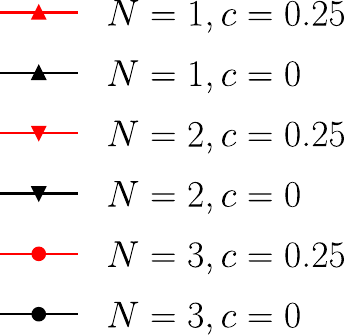}
	\end{minipage}
	\begin{minipage}[c]{0.85\textwidth}
		\centering{
			\subfloat{
				\begin{minipage}[b]{0.33\textwidth}
					\centering
					\includegraphics[width=1\linewidth]{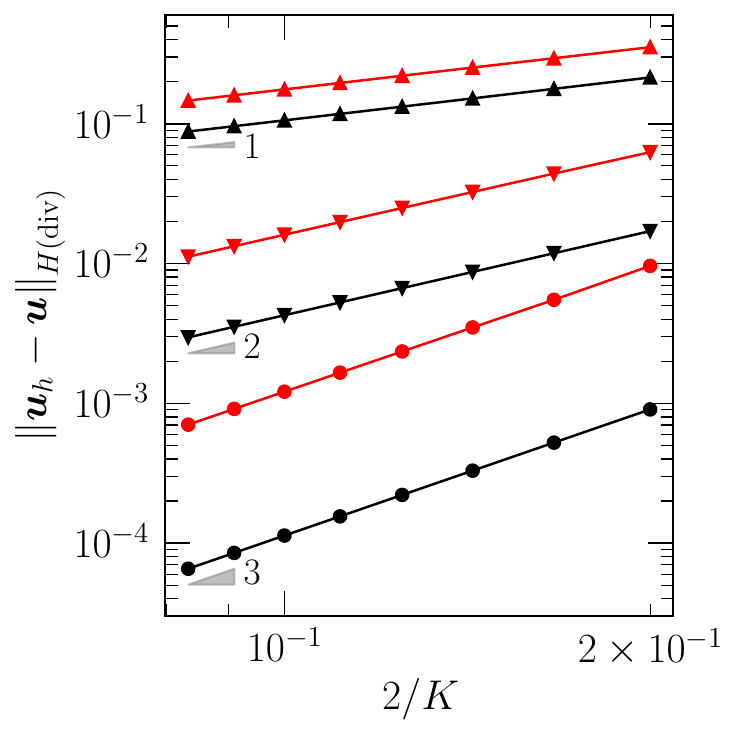}
				\end{minipage}
			}
			\subfloat{
				\begin{minipage}[b]{0.33\textwidth}
					\centering
					\includegraphics[width=1\linewidth]{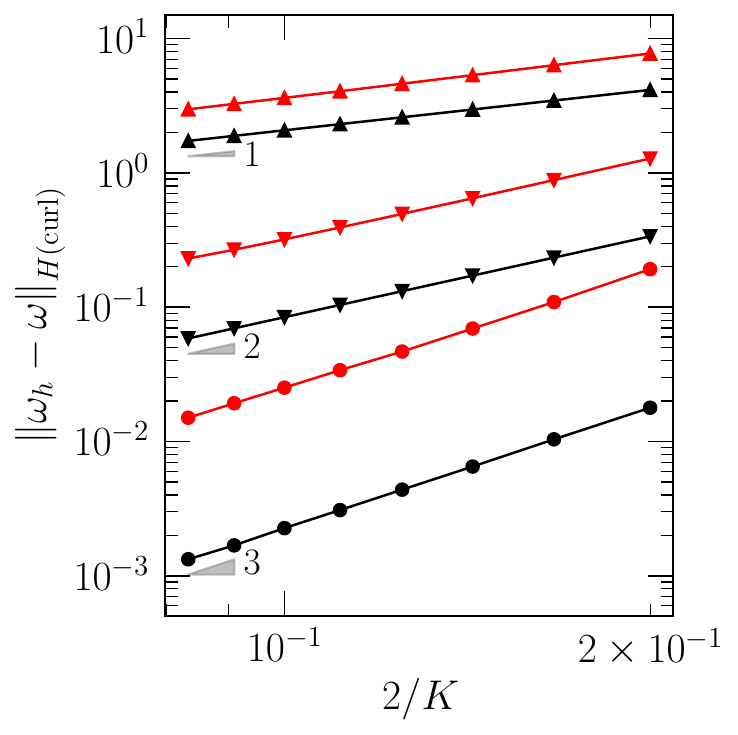}
				\end{minipage}
			}
			\subfloat{
				\begin{minipage}[b]{0.33\textwidth}
					\centering
					\includegraphics[width=1\linewidth]{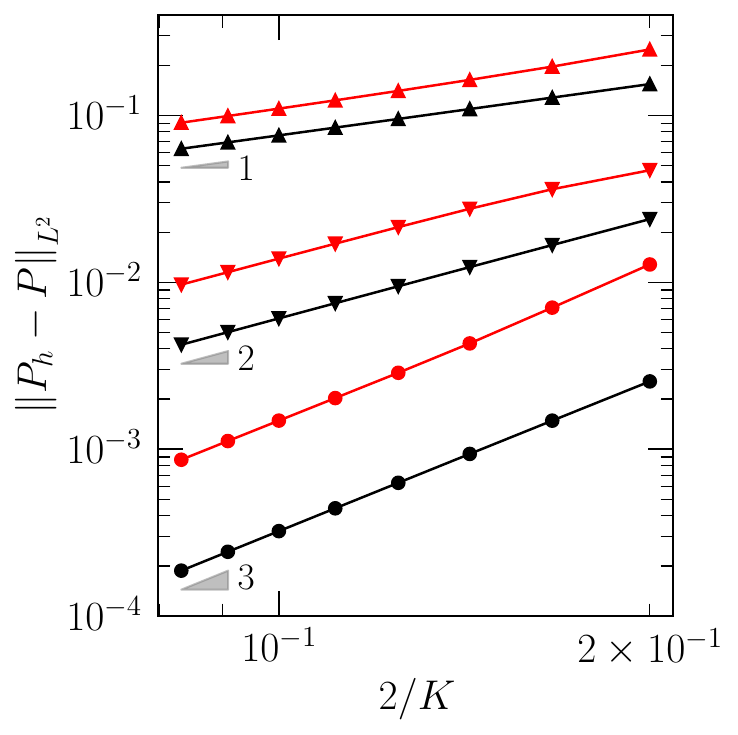}
				\end{minipage}
		}}
	\end{minipage}
	\caption{$ H(\mathrm{div}) $-error of $ \bu_{h} $, $ H(\mathrm{curl}) $-error of $ \omega_{h} $ and $ L^2 $-error of $ P_{h} $ at $ t=1 $ of the Taylor–Green vortex test under $ ph $-refinements for $ N\in\left\lbrace1,2,3\right\rbrace $, $c\in\left\lbrace0, 0.25\right\rbrace$, $ K\in\left\lbrace10,12,14\cdots,24\right\rbrace $, $\Delta t = \frac{1}{25}$ and $\Rn=100$.}
	\label{fig convergence tests 2d TGV}
\end{figure}

We test the accuracy of the method using a classic analytical solution of two-dimensional incompressible Navier-Stokes equations in the absence of external force, the Taylor–Green vortex, written as
\[
\begin{aligned}
	&u(x, y, t)=-\sin(\pi x)\cos(\pi y)e^{-2\pi^2t/\Rn } ,\\
	&v(x, y, t)=\cos(\pi x)\sin(\pi y)e^{-2\pi^2t/\Rn } ,\\
	&p(x, y, t) = \dfrac{1}{4} \left( \cos(2\pi x) + \cos(2\pi y)\right)e^{-4\pi^2t/\Rn} ,\\
	&\omega(x, y, t) = -2\pi\sin(\pi x)\sin(\pi y)e^{-2\pi^2t/\Rn }.
\end{aligned}
\]
The domain is set to $ \Omega = (x,y)\in[0,2]^2 $ with periodic boundary conditions. We use the meshes described by \eqref{eq:crazymesh} and solve the Taylor–Green vortex from $ t=t_{0}=0 $ to $ t=1 $ with $ \Rn=100 $ and $ \Delta t =\frac{1}{25} $. Results showing optimal convergence rates are presented in Fig.~\ref{fig convergence tests 2d TGV}.

\subsection{Conservation and dissipation tests: Shear layer roll-up} \label{subsec SLR}
\begin{figure}[h!]
	\centering{
		\subfloat[$ t=0 $]{
			\begin{minipage}[b]{0.3\textwidth}
				\centering
				\includegraphics[width=0.99\linewidth]{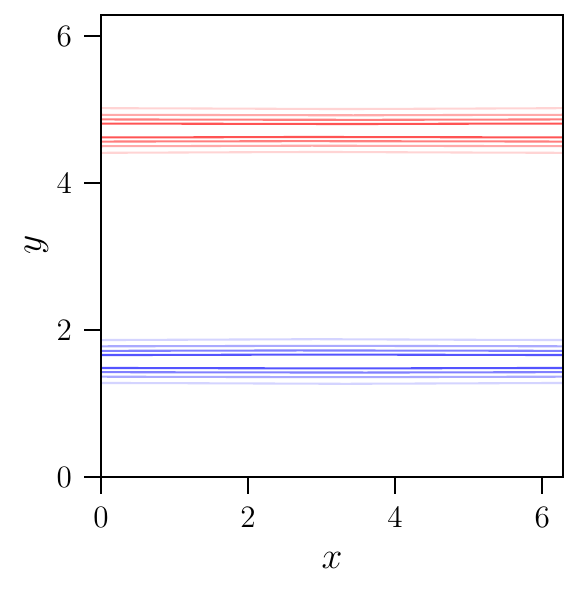}
			\end{minipage}
		}
		\subfloat[$ t=4 $]{
			\begin{minipage}[b]{0.3\textwidth}
				\centering
				\includegraphics[width=0.99\linewidth]{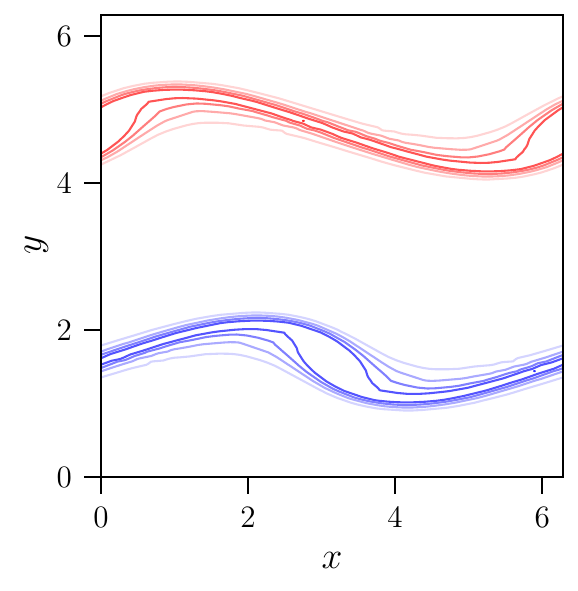}
			\end{minipage}
		}
		\subfloat[$ t=8 $]{
			\begin{minipage}[b]{0.3\textwidth}
				\centering
				\includegraphics[width=0.99\linewidth]{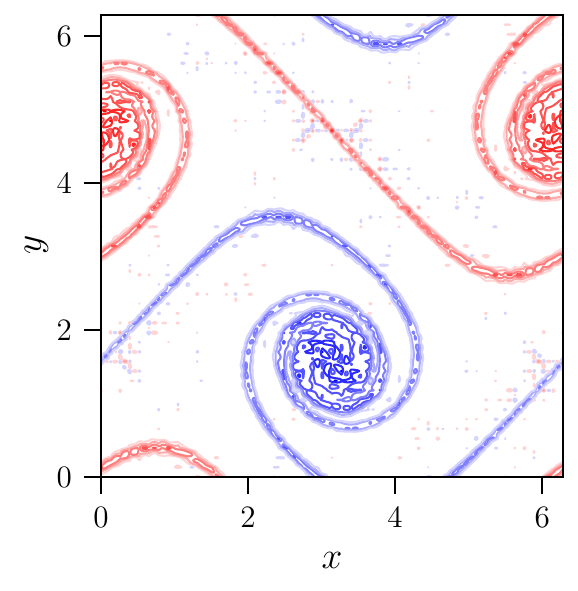}
			\end{minipage}
		}
		\caption{Vorticity field $ \omega_{h} $ of the ideal shear layer roll-up test at $ t\in\left\lbrace0,4,8\right\rbrace $ with contour lines for $ \omega_{h} \in \left\{\pm1, \pm2, \pm3, \cdots, \pm6\right\} $. The color scheme is from blue ($ -6\leftarrow\omega_{h} $) to red ($ \omega_{h}\rightarrow 6 $). The simulation is conducted for $N=2$, $c=0$, $K=48$, $\Delta t=\frac{1}{50}$.}
		\label{fig SLR vorticity field}
	}
\end{figure}

\begin{figure}[h!]
	\centering{
		\subfloat{
			\begin{minipage}[b]{0.5\textwidth}
				\centering
				\includegraphics[width=0.95\linewidth]{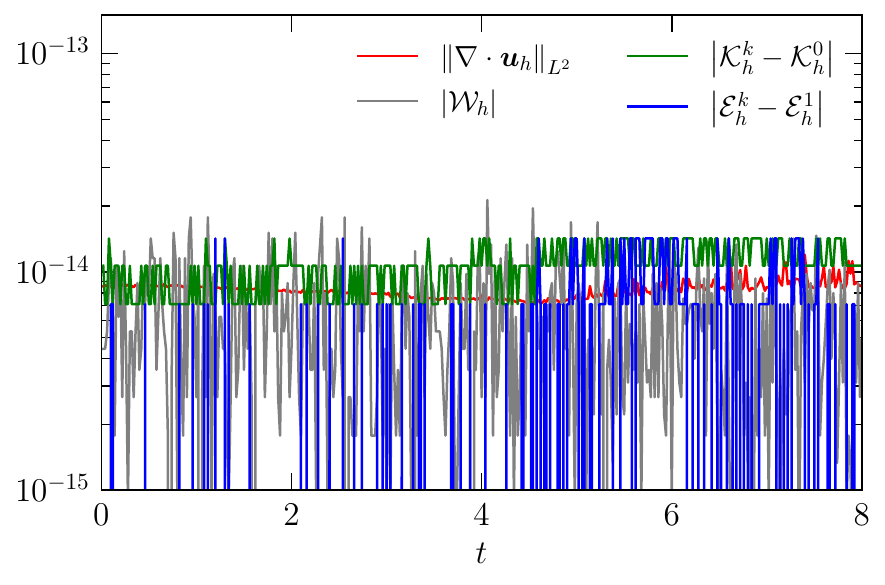}
			\end{minipage}
		}
		\subfloat{
			\begin{minipage}[b]{0.5\textwidth}
				\centering
				\includegraphics[width=0.95\linewidth]{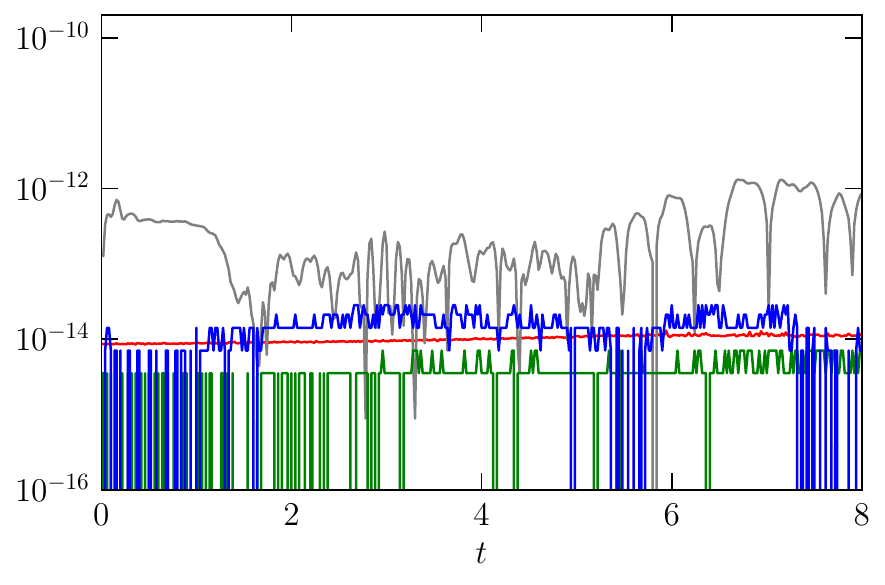}
			\end{minipage}
		}
		\caption{Discrete mass, energy, enstrophy and vorticity conservation over time of the ideal shear layer roll-up test for $N=2$, $c=0$ (left), $c=0.25$ (right), $K=48$ and $\Delta t=\frac{1}{50}$.}
		\label{fig SLR conservation}
	}
\end{figure}

The shear layer roll-up is a two-dimensional ideal flow whose initial condition is given by
\[
u^{0} = \left\lbrace
\begin{aligned}
	&\tanh\left(\dfrac{y - \pi/2}{\delta} \right), &\quad y \leq \pi \\
	&\tanh\left(\dfrac{3\pi/2 - y}{\delta}\right), &\quad y > \pi
\end{aligned}
\right.,
\qquad \qquad \qquad
v^{0} = \epsilon\sin(x),
\]
where $ \delta = \frac{\pi}{15} $ and $ \epsilon = 0.05 $, see \cite{Palha2017,sanderse2013energy}. The domain is $ \Omega=(x,y)\in[0, 2\pi]^2 $ with periodic boundary conditions. Meshes as described in \eqref{eq:crazymesh} for $c\in\left\{0,0.25\right\}$ and $K=48$ are used. The polynomial degree is set to $ N = 2 $ and the time interval is $ \Delta t = \frac{1}{50} $. The flow is computed from $ t=t_{0}=0 $ to $ t=8 $.
To limit the error caused by the Newton-Raphson method, the tolerance of outer iterations is set to $ 10^{-12} $. The vorticity field $ \omega_{h} $ at $ t\in\left\lbrace0,4,8\right\rbrace $ for $c=0$ is shown in Fig.~\ref{fig SLR vorticity field}. 
%\Yi{Note that in these tests, Gauss quadrature of a degree $N_{Q}=2N$ is used to compute discrete cross product terms.} 
In Fig.~\ref{fig SLR conservation}, results showing the conservation laws are satisfied to machine precision on both orthogonal and curvilinear meshes are presented. And throughout this section we use $ \left\|\nabla\cdot\bu_{h}\right\|_{L^{2}} $, i.e. the $ L^2 $-norm of $ \nabla\cdot\bu_{h} $, to identify mass conservation. Since the basis functions have normal continuity, if $ \left\|\nabla\cdot\bu_{h}\right\|_{L^{2}} =0 $ (to machine precision), pointwise mass conservation is satisfied everywhere.

\begin{figure}[h!]
	\centering{
		\subfloat{
			\begin{minipage}[b]{0.5\textwidth}
				\centering
				\includegraphics[width=0.95\linewidth]{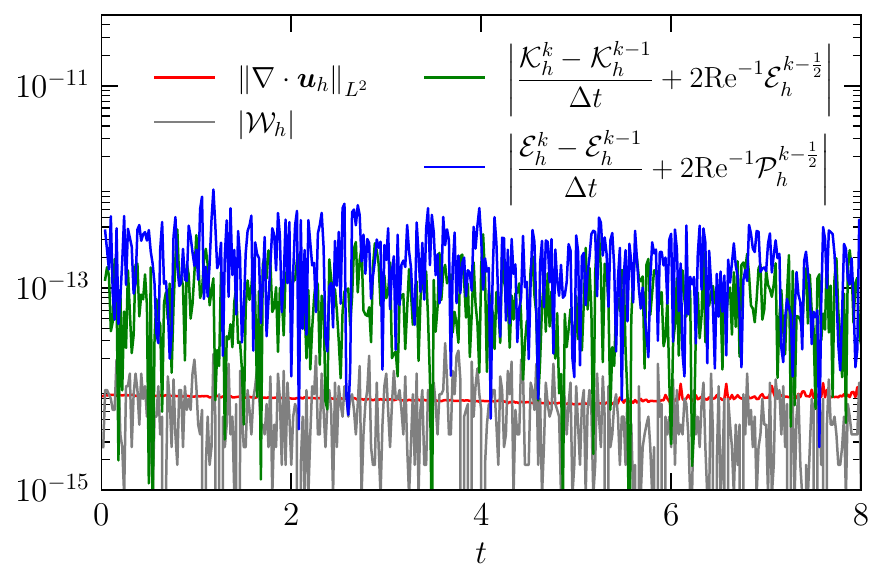}
			\end{minipage}
		}
		\subfloat{
			\begin{minipage}[b]{0.5\textwidth}
				\centering
				\includegraphics[width=0.95\linewidth]{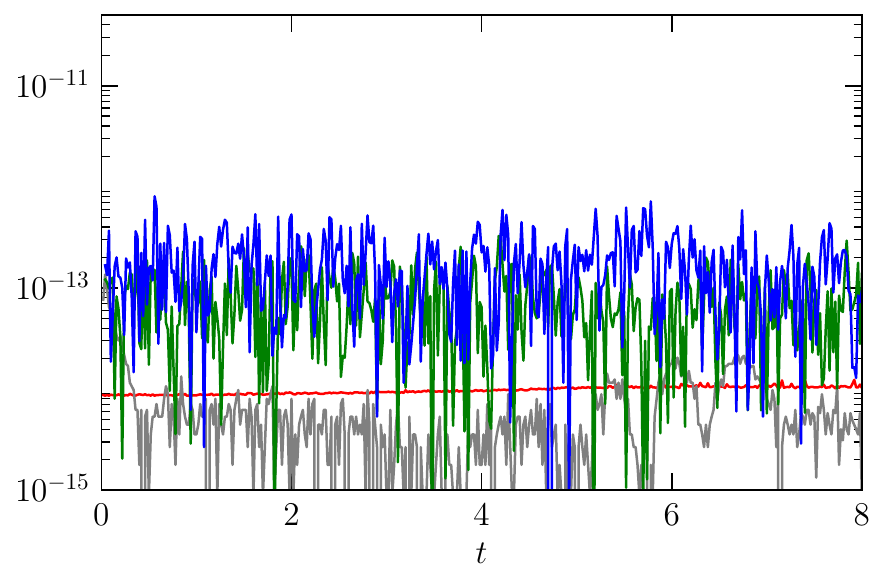}
			\end{minipage}
		}
		\caption{Discrete mass conservation, energy and enstrophy balances, and vorticity conservation over time of the viscous shear layer roll-up test for $N=2$, $c=0$ (left), $c=0.25$ (right), $K=48$, $\Delta t=\frac{1}{50}$ and $\Rn=500$.}
		\label{fig SLR dissipation}}
\end{figure}

We repeat the above test now for a viscous flow of $ \Rn = 500 $. The results are presented in Fig.~\ref{fig SLR dissipation}. It is seen that, to machine precision, (i) mass and vorticity conservation and (ii) energy and enstrophy balances, \eqref{Eq discrete energy dissipation} and \eqref{Eq discrete enstrophy dissipation}, are satisfied for both orthogonal and curvilinear meshes.

\subsection{No-slip boundary condition test: Normal dipole collision}
\label{subsec NDC}

\begin{figure}[h!]
	\centering
	\includegraphics[width=0.45\linewidth]{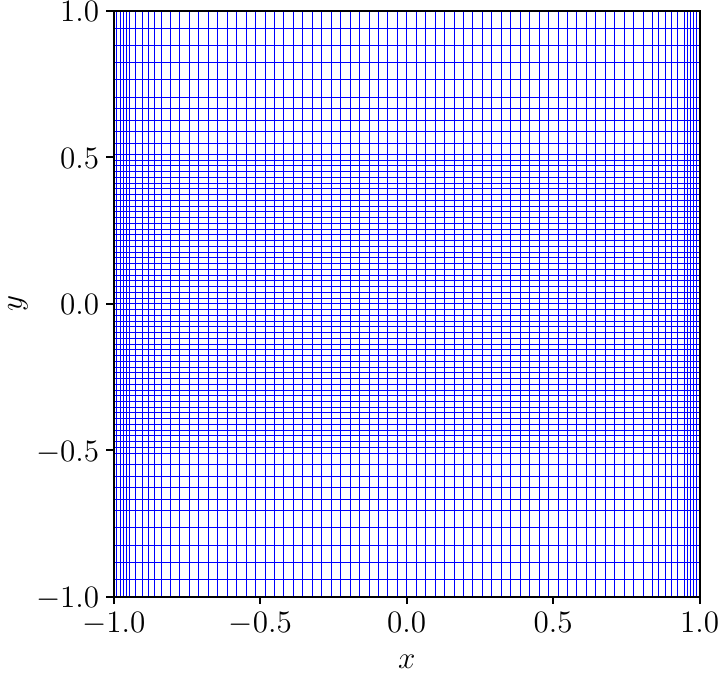}
	\caption{The mesh of $ 5148 $ elements used in the normal dipole collision test.}
	\label{fig NDC mesh}
\end{figure}

The normal dipole collision is a viscous flow in the domain $ \Omega=(x,y)\in[-1,1]^2 $ with no-slip boundary conditions on all four walls \cite{CLERCX2006245}. The unscaled initial velocity field, 
$ \boldsymbol{u}^{0}=\begin{bmatrix}
	u^{0}&v^{0}
\end{bmatrix} $, 
is given by
\[
\begin{aligned}
	&u^{0} = -\frac{1}{2}\left|\omega_{e}\right|\left( y-y_1\right) e^{-(r_{1}/r_{0})^2} +\frac{1}{2}\left|\omega_{e}\right|\left( y-y_2\right) e^{-(r_{2}/r_{0})^2},\\
	&v^{0} = -\frac{1}{2}\left|\omega_{e}\right|\left( x-x_2\right) e^{-(r_{2}/r_{0})^2}
	+\frac{1}{2}\left|\omega_{e}\right|\left( x-x_1\right) e^{-(r_{1}/r_{0})^2} ,
\end{aligned}
\]
where $ \left|\omega_{e}\right|=320 $, $ (x_{1}, y_{1})=(0, 0.1) $ and $ (x_{2}, y_{2})=(0, -0.1) $, $ r_{1} $ and $ r_{2} $ are distances to $ (x_{1}, y_{1}) $ and $ (x_{2}, y_{2}) $, respectively,  and $ r_{0}=0.1 $. This velocity field leads to a vorticity field expressed as
\[
\omega^{0} = \sum_{i\in\left\{1,2\right\}}\omega_{e,i}\left( 1-\left(\frac{r_{i}}{r_{0}}\right) ^2\right) e^{-(r_{i}/r_{0})^2},
\] 
where $ \omega_{e,1}=320$, $ \omega_{e,2}=-320 $, which is a combination of two monopoles centered at $ (x_{1}, y_{1}) $ and $ (x_{2}, y_{2}) $, respectively. The initial velocity is then scaled such that the initial kinetic energy is $ \mathcal{E}^{0} = 2 $. The scaling factor is $ f\approx0.936026 $. The corresponding initial enstrophy and palinstrophy are $ \mathcal{E}^{0}\approx 800 $ and $ \mathcal{P}^{0}\approx441855 $, respectively. For the present test, we use a non-uniform orthogonal mesh of $ 5148 $ elements, see Fig.~\ref{fig NDC mesh}, $ \Rn=625 $, polynomial degree $ N=2 $ and $ \Delta t = \frac{1}{200} $. This setup is similar to that in \cite{DeDiego2019} except that \cite{DeDiego2019} uses an unstructured mesh and much smaller time steps, $\Delta t= \frac{1}{2000}$.

\begin{figure}[!htb]
	\centering
	\begin{minipage}[c]{.92\textwidth}
		\centering{
			\subfloat[$ t=0 $]{
				\begin{minipage}[b]{0.32\textwidth}
					\centering
					\includegraphics[width=1\linewidth]{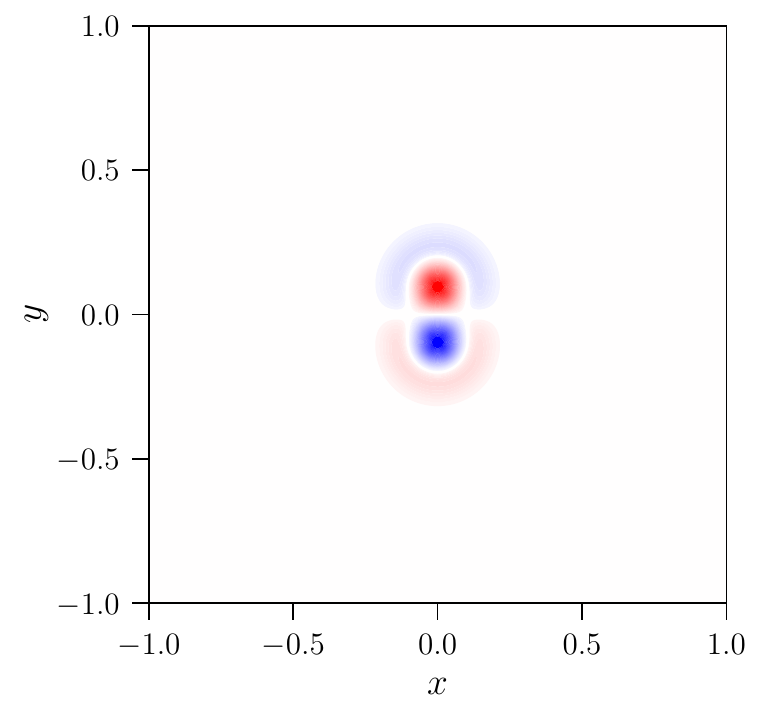}
				\end{minipage}
			}
			\subfloat[$ t=0.2 $]{
				\begin{minipage}[b]{0.32\textwidth}
					\centering
					\includegraphics[width=1\linewidth]{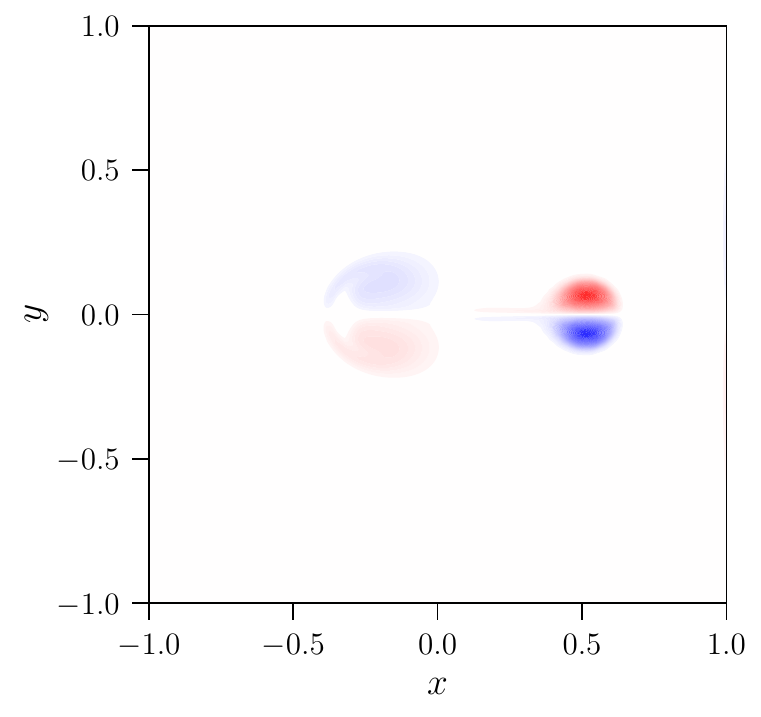}
				\end{minipage}
			}
			\subfloat[$ t=0.4 $]{
				\begin{minipage}[b]{0.32\textwidth}
					\centering
					\includegraphics[width=1\linewidth]{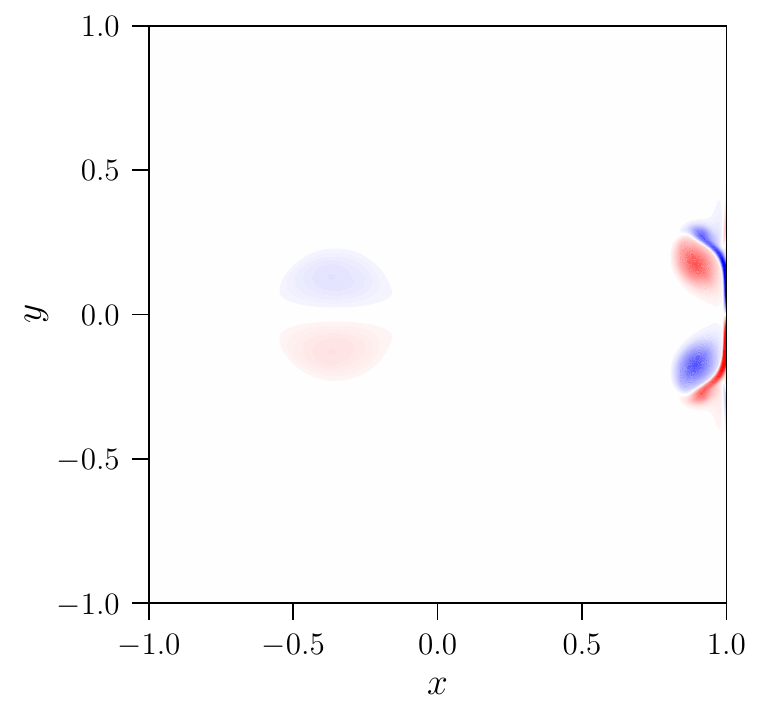}
				\end{minipage}
			}\\
			\subfloat[$ t=0.6 $]{
				\begin{minipage}[b]{0.32\textwidth}
					\centering
					\includegraphics[width=1\linewidth]{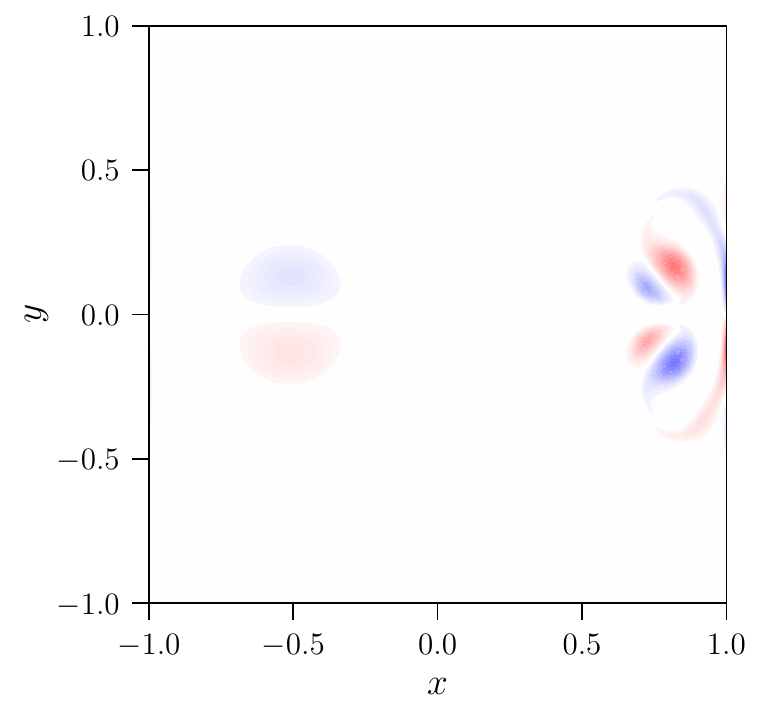}
				\end{minipage}
			}
			\subfloat[$ t=0.8 $]{
				\begin{minipage}[b]{0.32\textwidth}
					\centering
					\includegraphics[width=1\linewidth]{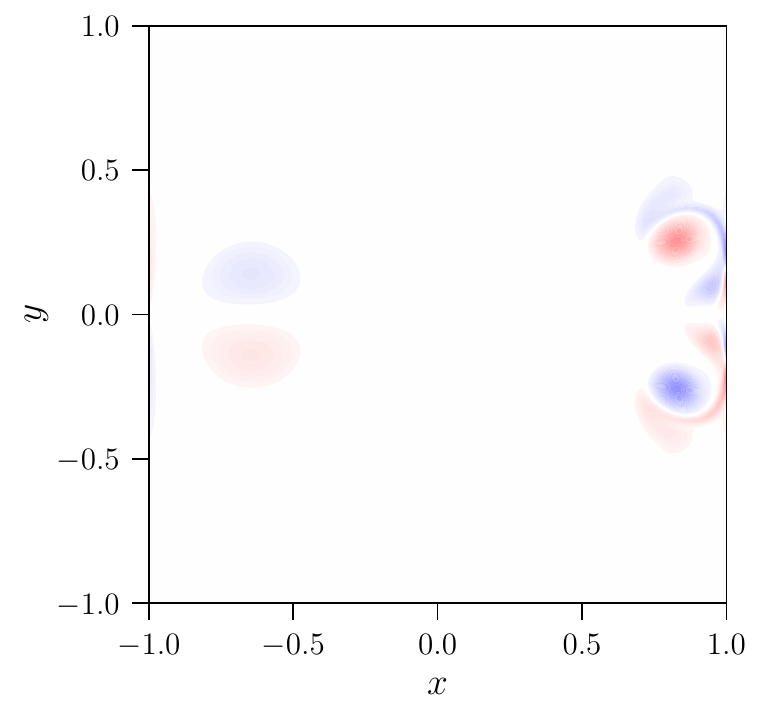}
				\end{minipage}
			}
			\subfloat[$ t=1 $]{
				\begin{minipage}[b]{0.32\textwidth}
					\centering
					\includegraphics[width=1\linewidth]{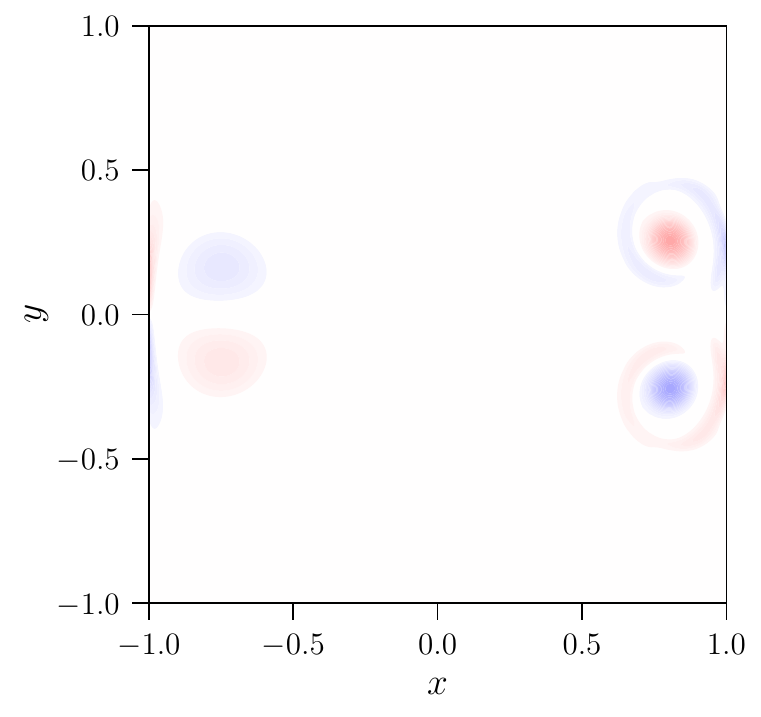}
				\end{minipage}
			}
		}
	\end{minipage}%
	\begin{minipage}[c]{0.08\textwidth}
		\centering
		\includegraphics[width=0.8\linewidth,]{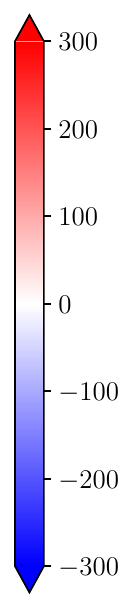}
	\end{minipage}
\caption{Vorticity field $ \omega_{h} $ at $ t\in\left\lbrace0, \frac{1}{5}, \frac{2}{5}, \cdots, 1\right\rbrace $ of the normal dipole collision test in a mesh of 5148 elements for $N=2$, $\Delta t = \frac{1}{200} $ and $\Rn=625$.}
\label{fig ndc mesh and vorticity}
\end{figure}

The dipole is initialized at $ t=t_{0}=0 $. It will move under a self-induced velocity in positive $x$-direction with an induced wake moving in the opposite direction. 
The simulation is performed until $ t=1 $. 
The vorticity field $ \omega_{h} $ at $ t\in\left\lbrace0, \frac{1}{5}, \frac{2}{5}, \cdots, 1\right\rbrace $ is presented in Fig.~\ref{fig ndc mesh and vorticity}. 

\begin{figure}[h!]
	\centering{
		\subfloat{
			\begin{minipage}[b]{0.33\textwidth}
				\centering
				\includegraphics[width=1\linewidth]{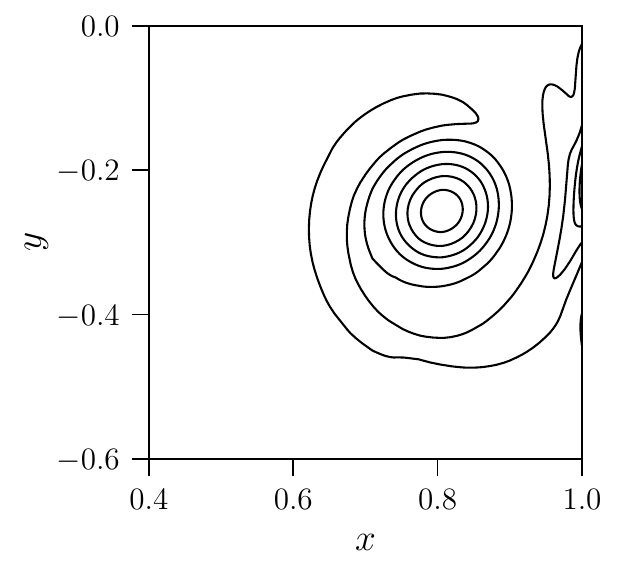}
			\end{minipage}
		}
		\subfloat{
			\begin{minipage}[b]{0.67\textwidth}
				\centering
				\includegraphics[width=0.94\linewidth]{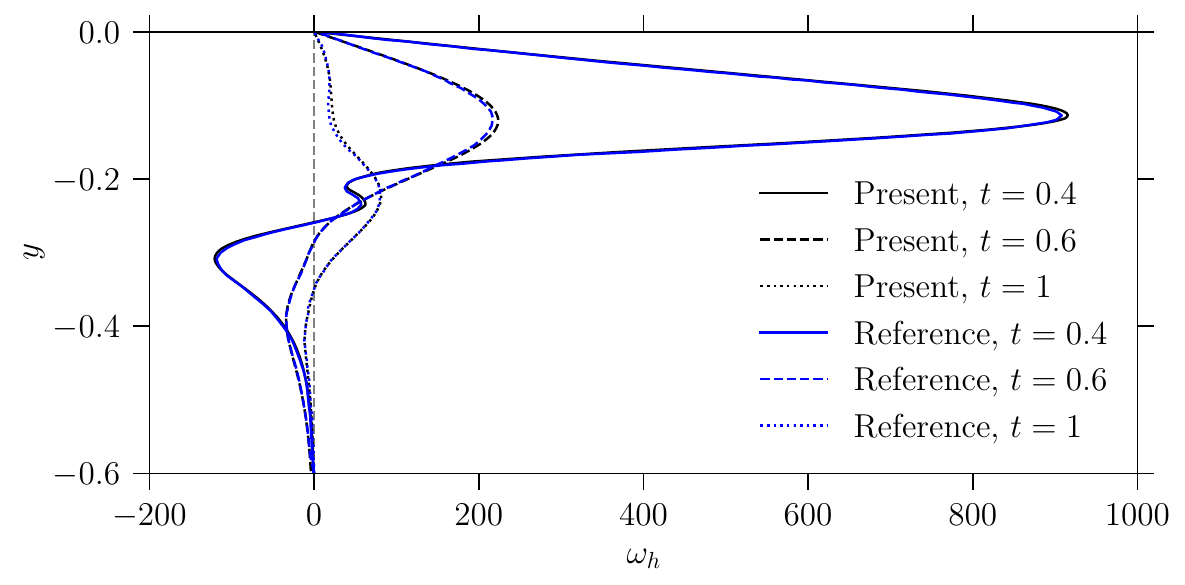}
			\end{minipage}
		}
		\caption{Vorticity field $ \omega_{h} $ in region $ \left(x,y\right)\in[0.4,1]\times[-0.6, 0] $ at $ t=1 $ with contour lines for $ \omega_{h}\in\left\lbrace -90, -70, -50, \cdots, 70\right\rbrace $ and on the boundary section $ (x,y)\in-1\times[-0.6, 0] $ at $ t\in\left\lbrace0.4,0.6,1\right\rbrace $ compared to reference results taken from \cite[Fig.~5]{CLERCX2006245} for $\Rn=625$. The present simulation has $145^2$ degrees of freedom for vorticity. The reference simulation uses a pseudospectral method and has $256^2$ degrees of freedom for vorticity.}
		\label{fig NDC wall vorticity}
	}
\end{figure}

\begin{figure}[h!]
	\centering{
		\subfloat{
			\begin{minipage}[b]{0.24\textwidth}
				\centering
				\includegraphics[width=1\linewidth]{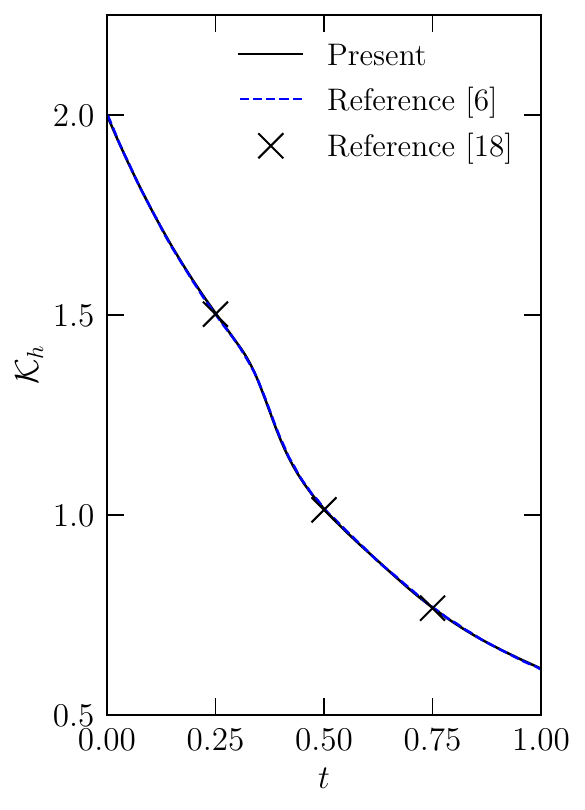}
			\end{minipage}
		}
		\subfloat{
			\begin{minipage}[b]{0.24\textwidth}
				\centering
				\includegraphics[width=1\linewidth]{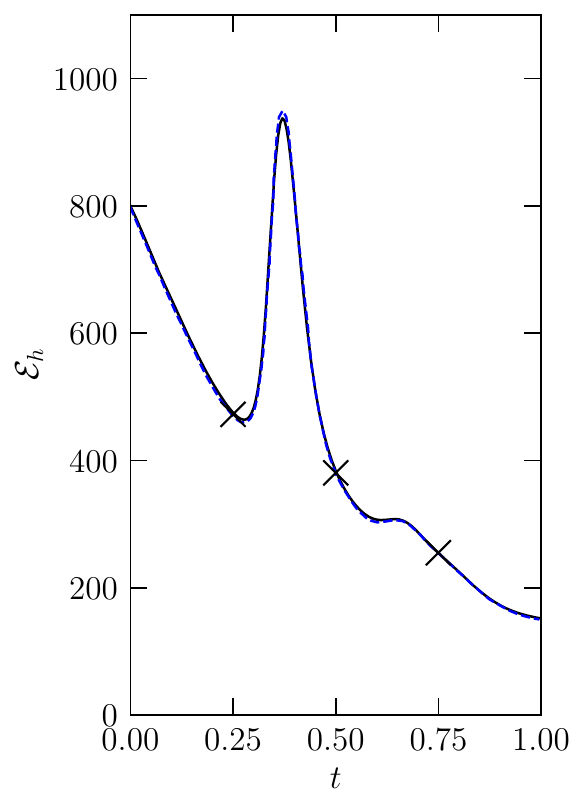}
			\end{minipage}
		}
		\subfloat{
			\begin{minipage}[b]{0.24\textwidth}
				\centering
				\includegraphics[width=1\linewidth]{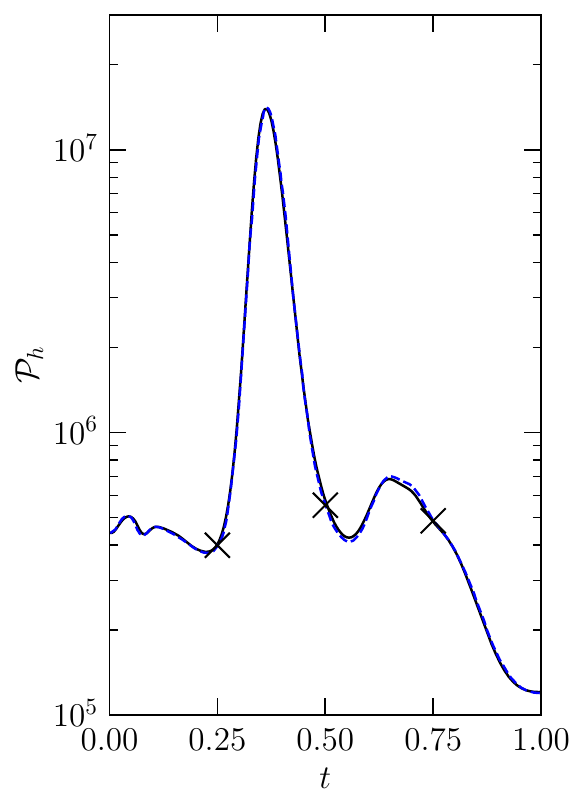}
			\end{minipage}
		}
		\subfloat{
			\begin{minipage}[b]{0.24\textwidth}
				\centering
				\includegraphics[width=1\linewidth]{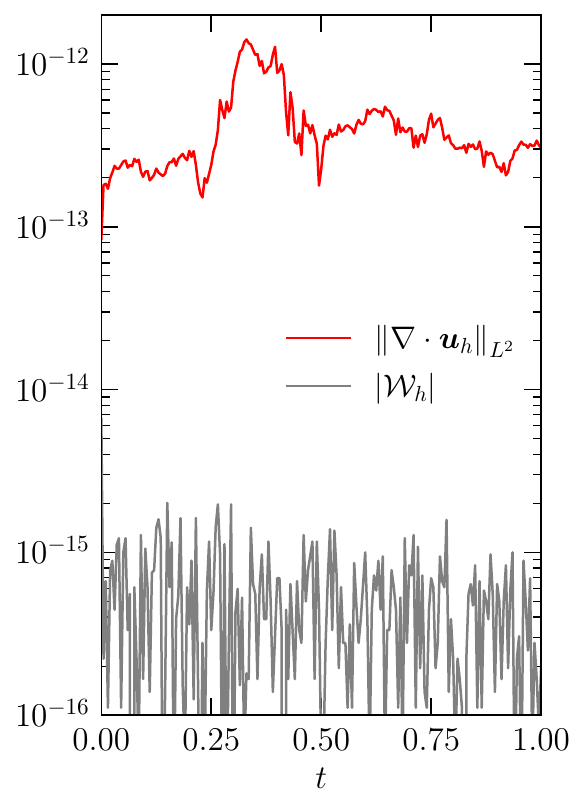}
			\end{minipage}
		}
		\caption{Discrete energy, enstrophy, palinstrophy over time compared to reference results 
			%using the original MEEVC method 
			taken from \cite[Fig.~8]{DeDiego2019} and reference results at $t\in\left\lbrace0.25,0.5,0.75\right\rbrace$ taken from \cite{CLERCX2006245}, and mass and vorticity conservation 
%			($\left\|\nabla\cdot\bu_{h}\right\|_{L^{2}} $) 
		over time.
%		The present simulation is conducted has $145^2$ degrees of freedom for vorticity. The reference simulation uses the original MEEVC scheme of $256^2$ degrees of freedom for vorticity.
	}
		\label{fig NDC dissipation and conservation}
	}
\end{figure}

The original MEEVC scheme uses indirect approaches to impose no-slip boundary conditions, and the suggested approach, called the kinematic Neumann approach, destroys vorticity conservation, see \cite[Fig.~9]{DeDiego2019}. In contrast, the present method can handle no-slip boundary conditions (and other general boundary conditions) naturally, see the formulation \eqref{Eq: form 1} or \eqref{Eq fully discrete}. In Fig.~\ref{fig NDC wall vorticity}, local distributions of vorticity field in region $ \left(x,y\right)\in[0.4,1]\times[-0.6, 0] $ at $ t=1 $ and on boundary section $ (x,y)\in-1\times[-0.6, 0] $ at different time instants (with comparisons to results in \cite{CLERCX2006245}) are shown.
The discrete energy, enstrophy and palinstrophy over time are presented and compared to results taken from \cite{DeDiego2019, CLERCX2006245} in Fig.~\ref{fig NDC dissipation and conservation} where mass and vorticity conservation is also shown. These results show an improved match with the reference than those in \cite[Fig.~10]{DeDiego2019} and also indicate that no-slip boundary conditions are correctly incorporated by the present method without destroying vorticity conservation.

\subsection{Convective term for enstrophy conservation} \label{Sub: ec}
This subsection provides numerical evidences for \eqref{Eq: key eq} as it is a key for the enstrophy balance of the present method. Given two random smooth scalar fields,
\[
\omega = 2\pi\sin(2\pi x + e)\sin(2\pi y + f)
\]
and
\[
\psi = 2 \pi \sin(2\pi x + g)\sin(2\pi y + h),
\]
where $e,f,g,h\in(0,1)$ are random real numbers, in the periodic unit square, $\Omega=(x,y)\in[0,1]^2$. In this domain, meshes as described in \eqref{eq:crazymesh} for $c\in\left\{0,0.25\right\}$ and $K=12$ are generated. $\omega$ and $\psi$ are projected to finite dimensional polynomial spaces as $\omega_{h}$ and $ \psi_{h} $, see \cite{Kreeft2011} or \cite{zhang2022phd} for details of these projections. The finite dimensional velocity is $\bu_{h} = \nabla\times \psi_{h}$ and thus $\nabla\cdot\bu_{h}=0 $. Then the trilinear form $a\left( \omega_{h},\boldsymbol{u}_{h},\nabla\times\omega_{h}\right)$ is computed with Gauss quadrature, see \cite{golub1969calculation}, of different degrees, $N_{Q}$. 

The results are presented in Table~\ref{Tab:1}. We can see that, for the orthogonal mesh ($c=0$), the trilinear form is zero (to machine precision) even when the quadrature is very inexact, for example, $N=4$ and $ N_{Q}=1$ (numerical quadrature of degree $N_Q$ is exact for polynomials of degree $N \leq 2N_Q - 1$). While for the curvilinear mesh $(c=0.25)$, as the metric of the mapping, see \eqref{eq:crazymesh}, cannot be captured by polynomials, the quadrature is always inexact and the trilinear form is still zero for quadrature degree that is significantly high. 
%As for the exceptions, see underlined values in Table~\ref{Tab:1}, although the authors do not understand, we suspect that it may related to the particular implementation we used and intend to study it carefully in the future, and obviously, they do not affect the conclusion that, 
These results support the statement that, for $\nabla\cdot\bu_{h}=0 $, the trilinear form  $a\left( \omega_{h},\boldsymbol{u}_{h},\nabla\times\omega_{h}\right)$ can be zero even with inexact numerical quadrature.

\begin{table}[h!]
	\centering
	\caption{$a\left( \omega_{h},\boldsymbol{u}_{h},\nabla\times\omega_{h}\right) $ for $c\in\left\{0,0.25\right\}$, $K=12$, $N\in\left\lbrace2,3,4\right\rbrace$ and $N_{Q}\in\left\lbrace1,2,3,4,5,6\right\rbrace$.}
	\label{Tab:1}
	\begin{tabular}{clrrrrrrr}
		\hline
		\multirow{2}{*}{$N_Q$} & \multicolumn{1}{c}{$c$} & \multicolumn{3}{c}{0}                                                 & \multicolumn{1}{l}{} & \multicolumn{3}{c}{0.25}                                              \bigstrut\\ \cline{3-5} \cline{7-9} 
		& \multicolumn{1}{c}{$N$} & \multicolumn{1}{c}{2} & \multicolumn{1}{c}{3} & \multicolumn{1}{c}{4} & \multicolumn{1}{c}{} & \multicolumn{1}{c}{2} & \multicolumn{1}{c}{3} & \multicolumn{1}{c}{4} \bigstrut\\ \hline
		1                      &                         & $1.05e-12$            & $4.72e-12$            & $-7.59e-11$           &                      & $\color{lightgray}-1.91e+02$           & $\color{lightgray}-6.93e+00$           & $\color{lightgray}1.25e+01$            \bigstrut\\
		2                      &                         & $-1.60e-12$           & $2.25e-12$            & $1.04e-12$            &                      & $-5.19e-12$           & $\color{lightgray}7.43e+00$            & $\color{lightgray}3.22e-01$            \bigstrut\\
		3                      &                         & $2.81e-13$            & $-1.01e-12$           & $3.65e-14$            &                      & $6.51e-12$            & $\color{lightgray}1.55e-03$            & $\color{lightgray}-7.90e-02$           \bigstrut\\
		4                      &                         & $1.65e-12$            & $7.96e-13$            & $1.01e-13$            &                      & $-7.46e-14$           & $-9.95e-14$           & $\color{lightgray}-1.48e-04$           \bigstrut\\
		5                      &                         & $-5.22e-12$           & $1.86e-12$            & $-7.21e-14$           &                      & $1.70e-12$            & $5.68e-13$            & $-1.42e-13$           \bigstrut\\
		6                      &                         & $3.46e-14$            & $5.59e-15$            & $1.71e-13$            &                      & $1.19e-13$            & $-2.52e-13$           & $2.23e-13$            \\ \hline
	\end{tabular}
\end{table}
\section{Conclusions} \label{Sec conclusions}
In this work, we present a  mass, energy, enstrophy and vorticity conserving (MEEVC) mixed finite element discretization for the rotational form of the incompressible Navier-Stokes equations on both orthogonal and curvilinear meshes. Comparing to the original MEEVC method, the present method uses a formulation of a single evolution equation and, more importantly, can naturally adapt no-slip boundary conditions without damaging vorticity conservation. However, it does not linearize the discrete systems as the original MEEVC scheme does; a more expensive nonlinear method has to be employed to solve the systems.

\section*{Acknowledgements}
The authors acknowledge dr. Andrea Brugnoli for helpful discussions.

%\section*{Appendix}

\bibliographystyle{elsarticle-num}
\bibliography{ref}

\end{document}